\documentclass[12pt,draftcls,onecolumn]{IEEEtran}

\ifCLASSINFOpdf
\else
\fi
\usepackage{graphicx}
\usepackage{epstopdf}
\usepackage{epsfig}
\usepackage{amsfonts}
\usepackage{amstext}
\usepackage{color}
\usepackage{empheq}
\usepackage{graphicx,colordvi}
\usepackage{amssymb,amsmath,amsopn,amsfonts,graphicx}
\usepackage{CJK}
\usepackage{float}

\usepackage{geometry}
\newtheorem{lemma}{\bf Lemma}[section]
\newtheorem{theorem}{\bf Theorem}[section]
\newtheorem{definition}{\bf Definition}[section]

\newtheorem{remark}{\bf Remark}[section]

\newtheorem{proposition}{\bf Proposition}[section]

\usepackage[justification=centering]{caption}

\begin{document}
%
%
%
%

\title{Linear Quadratic Mean Field Stackelberg Games: Open-loop and Feedback Solutions}

\author{
	Bing-Chang Wang, \emph{Senior Member, IEEE,}
	\thanks{This work was supported by the National Natural Science Foundation of China under Grants 61922051, 62122043 and 62192753.}
	\thanks{Bing-Chang Wang is with the School of Control Science and Engineering, Shandong University, Jinan,  China. (e-mail: bcwang@sdu.edu.cn) } Juanjuan Xu\thanks{Juanjuan Xu is with the School of Control Science and Engineering, Shandong University, Jinan, China. (e-mail: juanjuanxu@sdu.edu.cn) }, Huanshui Zhang, \emph{Senior Member, IEEE}
	\thanks{Huanshui Zhang is with the College of Electrical Engineering and Automation,
Shandong University of Science and Technology, Qingdao, China. (e-mail: hszhang@sdu.edu.cn) }
and Yong Liang
	\thanks{ Yong Liang is with the School of Information Science and Engineering, Shandong Normal University, Jinan, China. (e-mail: yongliang@sdnu.edu.cn) }
}

%
%

\markboth{Journal of \LaTeX\ Class Files}%
{Shell \MakeLowercase{\textit{et al.}}: Bare Demo of IEEEtran.cls for IEEE Journals}
%



\maketitle

\begin{abstract}
This paper investigates open-loop and feedback solutions of linear
quadratic mean field (MF) games with a leader and a large number
of followers. The leader first gives its strategy and then all the
followers cooperate to optimize the social cost as the sum of
their costs.  By variational analysis with MF approximations, we obtain a set of open-loop controls of players in terms of solutions to MF forward-backward stochastic differential equations (FBSDEs), which is further shown be to an asymptotic Stackelberg equilibrium.
By applying the matrix maximum principle,  a set of decentralized feedback
strategies is constructed for all the players.
For open-loop and feedback solutions,
 the corresponding optimal costs of all players
 are explicitly given by virtue of the solutions to two Riccati equations, respectively.
 The performances of two solutions are compared by the numerical simulation.

\end{abstract}
%
%
%
%
%
%
%
%
%
\section{Introduction}
\subsection{Background and Motivation}
Mean field (MF) game theory is proposed to effectively design decentralized strategies for large-population models where the effect of each player is negligible while the influence of the whole population is significant \cite{LL07, hcm2007a, wbvr2008}.
The main methodology of MF games is to replace the interactions among all players by population aggregation effect, which structurally depicts the MF interactions in large population systems. 
Along this line, decentralized solutions may be obtained accordingly \cite{BFY13}, \cite{C14}, \cite{CD18}. MF games have been found wide applications in many areas such as dynamic production adjustment \cite{bibitem181}, vaccination games \cite{bibitem5}, and charging control of electric vehicles \cite{MCH13, TKG20},  resource allocation in internet of things \cite{bibitem7}, and etc.

We now present a compact review for the studies of MF games. 
First, depending on the state and cost setup of an MF game, it can be classified into the linear-quadratic (LQ) type or the nonlinear type. The LQ type has weak coupling through population state average, and it is commonly adopted in MF studies because of its analytical tractability and close connection to practical applications. Some relevant works include \cite{hcm2007a, lz2008, bj2012a, bp2014, BSY16, HH16}. 
Meanwhile, the nonlinear type of MF game 
is also of great importance because of its modeling generality and close relation to physics. A large body of work is devoted to this type; see \cite{LL07, CD13, L16}.
Second, depending on their system hierarchy, MF games can be classified into homogeneous, heterogeneous, or mixed game.
In a homogeneous game, all agents are symmetric 
and minor. In a heterogeneous game, all agents are still minor but may demonstrate diversity in their system coefficients.
A mixed game is
more distinctive especially in its decision hierarchical structure, rather than merely
in the coefficient datum. It involves some major agents by imposing a dominant impact
on all minor agents. 
A mixed game is a realistic setup for modeling a monopoly in economic dynamics. See \cite{H10, NH12, WZ12} for mixed games. Hence, it has attracted considerable research attentions.
For instance, \cite{H10, NH12} investigated LQ mixed MF games with a major agent and various minor agents, and provide $\epsilon$-Nash equilibrium strategies. \cite{WZ12} study a mixed game in a discrete-time case.
For nonlinear MF mixed games, see \cite{NC13, BLP13, SC16}. For more comprehensive literature, readers may refer to \cite{C12}
for an overview of MF games, 
and the surveys \cite{GS13}, \cite{C14}. 
We also draw attention to the recent monographs such as \cite{BFY13, CD18}.

As a typical class of mixed games, 
the Stackelberg 
game contains at least two hierarchies of players. One hierarchy of players are defined as leaders with a dominant position and the other players are defined as followers with a subordinate position. The leader has the priority to announce a strategy first and then the followers seek strategies to minimize their costs with response to leader's strategy. Taking account of followers' optimal responses, the leader chooses a strategy to optimize its cost. 
The study of Stakelberg games has a long history, and most early works focused on Stakelberg models with a leader and a follower (see e.g., \cite{sc1973}, \cite{y02}, \cite{bbs2010}, \cite{XZC15}, \cite{LZ19}, \cite{HWX21}, \cite{SXZ22}). Recently, MF Stackelberg games have attracted  much research interests \cite{bcy2015a, BCC20, bj2014a}, \cite{mb2018},  \cite{YH21}. Bensoussan et al. \cite{bcy2015a} studied MF Stackelberg games with delayed instructions, and further generalized to the case with the heterogeneous delay effect from leader's
action \cite{BCC20}.  Moon and Basar \cite{MB15, mb2018} investigated continuous-time MF-LQ Stackelberg games by the fixed-point method. Authors in \cite{BAA20} gave the saddle-point strategy of the
 minmax MF team
for a leader-follower network.
Besides, the work \cite{bj2014a} considered discrete-time hierarchical MF games with
tracking-type costs and gave the feedback $\varepsilon$-Stackelberg equilibrium. Authors in \cite{YH21} 
investigated feedback strategies of MF-LQ Stackelberg games by solving the master equations of the  limit
model.

Different from noncooperative games, the social optimization problem is a joint decision problem where all the players have the same goal and work cooperatively to optimize the social cost. This is regarded as a type of team decision problem \cite{hch1972}. Huang \emph{et al.} \cite{hcm2012} considered social optima in MF-LQ control, and provided an asymptotic team-optimal solution. The work \cite{bj2017b} investigated the MF social optimal problem where the Markov jump parameter appears as a common source of randomness. 
For further literature, see \cite{hn2016} for social optima in mixed games, \cite{AM15} for team-optimal control with finite population and partial information, \cite{SY20} for stochastic dynamic teams and their mean-field limit.
 \subsection{Contribution and Novelty}
This paper investigates MF-LQ Stackelberg games with a leader and many followers. The leader first gives his strategy and then all the followers cooperate to optimize the \emph{social cost}--the sum of their individual costs. For instance, consider an example of macroeconomic regulation, where the government is the leader, and companies in a market are followers \cite{P77}.  Another example is decentralized hierarchical planning of plugged-in electric vehicles (PEVs), which is modelled as
an MF reverse Stackelberg game in \cite{TKG20}.
Different from \cite{BCC20}, \cite{bcy2015a}, \cite{mb2018}, \cite{bj2014a}, our model involves both the leader's \emph{control} and population state average 
in dynamics and costs of followers. This implies that the leader and followers are highly interactive and coupled.
Until now, most relevant works studied open-loop control of  MF Stackelberg games, and rare
works focused on feedback strategies. Particularly, the relationship between open-loop and
feedback strategies is still unclear.



In this paper, we systematically study open-loop and feedback solutions to MF Stackelberg games by decoupling
forward-backward stochastic differential equations (FBSDEs).
We first consider the open-loop solution of MF Stackelberg games. Under a given control of the leader, we solve a centralized social control problem by the variational analysis, which leads to a system of high-dimensional FBSDEs. By MF approximation, we obtain a set of (strict) open-loop controls for followers with the help of an MF FBSDE, which can be implemented offline.
After applying the followers' strategies,
By solving an optimal control problem driven by  FBSDEs,
we  obtain a decentralized strategy of the leader.
By perturbation analysis, the proposed decentralized strategy is shown to be an $(\varepsilon_1,\varepsilon_2)$-Stackelberg equilibrium. Furthermore, we obtain the asymptotic optimal costs of players in terms of the solutions to Riccati equations. Next, we study the feedback solution of MF Stackelberg games. Different from the open-loop control, we  presume that the strategy of the leader has a feedback form. Fixing the feedback gain of the leader, we first solve a
centralized social control problem for followers by decoupling high-dimensional FBSDEs. By applying the matrix maximum principle with MF
approximations, we solve the optimal control problem
for the leader and construct decentralized feedback
strategies for all players. By the technique of  completing the square, we show that the proposed decentralized
strategy is a feedback $(\varepsilon_1,\varepsilon_2)$-Stackelberg equilibrium and further give the explicit form of the corresponding costs of players.
Finally, the performances under two solutions are compared by the numerical simulation. 

The main contributions of the paper are summarized as follows.
\begin{itemize}
	\item By applying variational analysis and MF approximations, we obtain a set of (strict) open-loop control laws in terms of solutions to MF FBSDEs, which is further shown be to an asymptotic Stackelberg equilibrium.

	\item  With the leader's feedback gain fixed, we first obtain centralized strategies of followers by decoupling high-dimensional FBSDEs. By applying the matrix maximum principle,  a set of decentralized feedback
	strategies is constructed for all players.

	\item
	For open-loop and feedback Stackelberg solutions,
	the corresponding  costs of players
	are explicitly given in terms of the solutions to Riccati equations, respectively.
	By the numerical simulation,
	it is shown that the performances of two solutions are different
	for the leader and  followers.
	
\end{itemize}

Below we highlight some key differences between this work and previous works.

(i) Different from \cite{bcy2015a, mb2018, BCC20}, we consider the case that followers cooperate to optimize the social cost where the dynamics of each player involves the state averages of followers. For social control problems, additional terms  need to be introduced when the auxiliary control problem is obtained by the fixed-point approach \cite{hcm2012}, \cite{bj2017b}, \cite{HWX21}. This is quite  different from  the case of noncoopratative games.
In this paper, we adopt the direct approach \cite{HZ20, WZZ20, LL07} to design decentralized strategies in terms of the solution to  8$n$-dimensional FBSDEs. Compared with previous work \cite{HWX21} resorting to tackling 10$n$-dimensional FBSDEs,  
less
computation is required. 


(ii) There are rare works on feedback solutions to MF  Stackelberg games. The work \cite{YH21}  treated the MF term as a part of agents' dynamics and employed dynamic programming to design feedback strategies. However, their scheme will be very complicated to solve the social control problem of followers.   
By applying the direct approach and the matrix
maximum principle, we construct a set of decentralized feedback strategies for all players. Different from previous works,
a cross term is introduced in derivation due to the appearance of MF effect.


(iii) To our best knowledge,   
the relationship between open-loop and feedback strategies is still unclear for  MF  Stackelberg games. Under 
different $\varepsilon$-Stackelberg equilibria,
the corresponding costs of players
are explicitly given by virtue of the solutions to Riccati equations, respectively.
By the numerical simulation, it is shown that the feedback solution generally outperforms the open-loop solution for followers, 
while the opposite is true for the leader.

\subsection{Organization and Notation}
{The paper is organized as follows. In Section II, we formulate the problem of LQ-MF
	 Stakelberg games. 
	In Section III,  we first obtain a set of open-loop control laws in terms of  MF FBSDEs, and give its feedback representation by virtue of Riccati equations.
	In Section IV, we design the feedback strategies of MF Stakelberg games and give the corresponding costs of all players. In Section V, we provide a numerical example to compare the performance of two solutions. 
	Section VI concludes the paper. }

\emph{Notation}: The following notation will be used throughout this paper. $\|\cdot\|$
denotes the Euclidean vector norm or matrix spectral norm. For a vector $z$ and a matrix $Q$, $\|z\|_Q^2= z^TQz$; $Q>0$ ($Q\geq0$) means that the matrix $Q$ is positive definite (positive semi-definite).
For two vectors $x,y$, $\langle x,y\rangle=x^Ty$.
$L^2
_{\cal F}(0, T; \mathbb{R}^k)$ is the space of all $\mathcal{F}_t$-adapted $\mathbb{R}^k$-valued processes $x(\cdot)$ such that
$\mathbb{E}\int_0^T\|x(t)\|^2dt<\infty$.
For convenience of  presentation, we use $C, C_1,C_2,\cdots$ to
denote generic positive constants, which may vary from place to place.

\section{Problem Formulation}

Consider a large-population system with a leader and $N$ followers. The states of the leader $\mathcal{A}_0$ and the $i$th follower $\mathcal{A}_i$, $1\leq i\leq N$ evolve by the following linear SDEs: 
\begin{equation}\label{eq1}
\left\{
\begin{aligned}
dx_0(t)=&[A_0x_0(t)+B_0u_0(t)
+G_0x^{(N)}(t)]dt
+D_0dW_0(t),\\
dx_i(t)=&[Ax_i(t)+Bu_i(t)+Gx^{(N)}(t)+Fx_0(t)+B_1u_0(t)]dt
+DdW_i(t),\\
x_0(0)=&\ \xi_0, \quad \quad x_i(0)=\xi_{i}, \quad i=1,2,\cdots,N,
\end{aligned}
\right.
\end{equation}\\
where $x_i\in \mathbb{R}^n,u_i\in \mathbb{R}^m$ are the state and control of agent $i$, $i=0,1,\cdots,N$, respectively.
$x^{(N)}(t)\triangleq\frac{1}{N}\sum_{i=1}^{N}x_i(t)$ is the state average of all the followers. $W_i(\cdot), i=0,\cdots,N$ are a sequence of $d$-dimensional standard Brownian motion defined on $(\Omega,\mathcal{F},\{\mathcal{F}_t\}_{0\leq t\leq T}, \mathbb{P})$. 
 Let $\mathcal{F}_t=\sigma(\xi_0, \xi_i, W_0(s), W_i(s), \\ s\leq t, 1\leq i\leq N)$). Denote $\mathcal{F}_t^0=\sigma(\xi_0, W_0(s), s\leq t)$ and $\mathcal{F}_t^i=\sigma(\xi_0, W_0(s), \xi_i, W_i(s), 0\leq s\leq t)$ for $i=1,\cdots, N$. 
The \emph{decentralized} 
admissible control set for the leader is defined by
$
\mathcal{U}_{d}^0=\big\{u_0|u_0(t)
\in L_{\mathcal{F}_t^0}^2(0,T;\mathbb{R}^m)\big\}
$.
The 
 \emph{decentralized} admissible control set for all the followers are defined by
\begin{equation}\nonumber
\begin{split}
\mathcal{U}_{d}=&\Big\{(u_1,\cdots,u_N)|u_i(t) \in L_{\mathcal{F}_t^i}^2(0,T;\mathbb{R}^m), 1\leq i \leq N\Big\}.
\end{split}
\end{equation}
Also, we define the centralized control sets for the leader and followers as $
\mathcal{U}_{c}^0=\big\{u_0|u_0(t)
\in L_{\mathcal{F}_t}^2(0,T;\mathbb{R}^m)\big\}
$ and
\begin{equation*}
\begin{split}
  \mathcal{U}_c =& \Big\{(u_1,\cdots,u_N)|u_i(t)\in L_{\mathcal{F}_t}^2(0,T;\mathbb{R}^m), 1\leq i \leq N\Big\}.
\end{split}
\end{equation*}

For the leader $\mathcal{A}_0$, the cost functional is defined by
\begin{equation}\label{eq2}
\begin{aligned}
{J}_0(u_0,u)=&\mathbb{E}\int_{0}^{T}\big[\|x_0(t)- \Gamma_0x^{(N)}(t)\|_{Q_0}^{2}  +\|u_0(t)\|_{R_0}^{2}\big]dt
+\mathbb{E}\big[\|x_0(T)-\bar{\Gamma}_0x^{(N)}(T)\|^2_{H_0}\big],
\end{aligned}
\end{equation}
where $u=(u_1,\cdots,u_N)$; $Q_0$ and $R_0$ are constant matrices with proper dimensions. For the $i$th follower  $\mathcal{A}_i$, the cost functional is defined by
\begin{equation}\label{eq3}
\begin{aligned}
J_i(u_0,u)=&\mathbb{E}\int_{0}^{T}\big[\|x_i(t)- \Gamma x^{(N)}(t)-\Gamma_1x_0(t)\|_{Q}^{2}+\|u_i(t)\|_{R}^{2}+2u_i^T(t)Lu_0(t)\cr
&+\|u_0(t)\|_{R_1}^{2}\big]dt
+\mathbb{E}\big[\|x_i(T)-\bar{\Gamma}x^{(N)}(T)-\bar{\Gamma}_1x_0(T)\|^2_{H}\big],
\end{aligned}
\end{equation}
where $Q, R,L$ and $R_1$ are constant matrices with proper dimensions. All the followers cooperate to minimize their social cost functional, which is denoted by
\begin{equation}\label{eq4}
{J}_{\rm soc}^{(N)}(u_0,u)=\sum_{i=1}^{N}{J}_i(u_0,u).
\end{equation}

Now we introduce the following assumptions.

(\textbf{A1}) $x_i(0),i=0,1,\cdots,N$ are a sequence of independent random variables.  
$\mathbb{E}x_i(0)=\bar{\xi}$, $1\leq i\leq N$. There exists a constant $c_0$ such that $\sup_{1\leq i \leq N}\mathbb{E}\|x_i(0)\|^2\leq c_0$. Furthermore, $\{x_i(0)\}$ and $W_i(t),i=0,1,2,\cdots,N$ are independent with each other. 

(\textbf{A2}) $Q_0\geq 0$, $H_0\geq 0$, $R_0>0$ and $Q\geq 0$, $H\geq 0$, $R>0$.

In this paper, we discuss the Stackelberg 
solution of the leader-follower  game. The leader, $\mathcal{A}_0$, holds
a dominant position in the sense that it first
announces his strategy $u_0$, and enforces on $\mathcal{A}_i$, $1 \leq i \leq N$.
The $N$ followers then respond by cooperating to optimize their social cost (\ref{eq4})
 under the leader's strategy.
 In
this process, the leader takes into account the rational reactions
of followers.  Note that in this framework,
the followers know the leader's strategy and state.

 Due to accessible information restriction and high computational
complexity, one generally cannot attain centralized Stackelberg equilibrium, but asymptotic Stackelberg equilibrium under decentralized information structure. We now introduce the definition of the open-loop ($\epsilon_1,\epsilon_2$)-Stackelberg equilibrium.

 \begin{definition}
  A set of control laws $(u_0^*,
	u_1^*,\cdots,u_N^*)$
	is an open-loop ($\epsilon_1,\epsilon_2$)-Stackelberg equilibrium 
	if we have:
	
	(i) For a given strategy $u_0^*\in \mathcal{U}^0_d$, $t\in [0,T]$, $u^*=(u_1^*,\cdots,u^*_N)$ is an $\epsilon_1$-optimal response if $u^*$ has $\epsilon_1$-team optimality, i.e.,
$$\frac{1}{N}J^{(N)}_{\rm soc}(u_0^*,u^*)\leq\frac{1}{N}J^{(N)}_{\rm soc}(u_0^*,u)+\epsilon_1, \ \hbox{for any}\ u\in \mathcal{U}_c, $$

	(ii) For any $u_0\in \mathcal{U}_d^0$, $J_0(u^*_0,u^*)\leq J_0(u_0,u)+\epsilon_2$, where $u^*,u$ are $\epsilon_1$-optimal responses to strategies $u_0^*,u_0$, respectively.


 \end{definition}

Motivated by \cite{mb2018, HZ20,WZZ20}, we consider feedback strategies with the following form:
\begin{equation}\label{eq3b}
  \left\{
\begin{aligned}
u_0(t)=&P_0(t)x_0(t)+\bar{P}(t)\bar{x}(t),\cr
  u_i(t)=&K(t)x_i(t)+\bar{K}(t)\bar{x}(t)+K_0(t)x_0(t), \ i=1,\cdots,N
\end{aligned}
\right.
\end{equation}
where $P_0,\bar{P},{K},\bar{K},{K}_0\in L_2(0,T;\mathbb{R}^{n\times n})$; $x_0,x_i$ and $\bar{x}$ satisfy
\begin{equation}\label{eq3c}
\left\{
\begin{aligned}
dx_0(t)=&[A_0x_0(t)+B_0(P_0(t)x_0(t)+\bar{P}(t)\bar{x}(t))
+G_0x^{(N)}(t)]dt
+D_0dW_0(t),\\
dx_i(t)=&[Ax_i(t)+B(K(t)x_i(t)+\bar{K}(t)\bar{x}(t)+K_0(t)x_0(t))+Gx^{(N)}(t)+Fx_0(t)\cr
&+B_1(P_0(t)x_0(t)+\bar{P}(t)\bar{x}(t))]dt
+DdW_i(t),\cr
d\bar{x}(t)=&\big\{[A+G+B(K(t)+\bar{K}(t))+B_1\bar{P}(t)]\bar{x}(t)+[F+BK_0(t)+B_1P_0(t)]x_0(t)\big\}dt,\cr
x_0(0)=&\xi_0, \quad x_i(0)=\xi_{i}, \quad i=1,2,\cdots,N, \quad  \bar{x}(0)=\bar{\xi}.
\end{aligned}
\right.
\end{equation}
In the above, $\bar{x}$ is an approximation of $x^{(N)}$ for $N\to \infty$.  We now introduce the definition of the feedback ($\epsilon_1,\epsilon_2$)-Stackelberg equilibrium.

 \begin{definition}
 	A set of control laws $(\hat{u}_0,
	\hat{u}_1,\cdots,\hat{u}_N)$
	is a feedback ($\epsilon_1,\epsilon_2$)-Stackelberg equilibrium 
	if the following hold:

	(i)  When the leader announces a strategy $\hat{u}_0=P_0x_0+\bar{P}\bar{x}$ at time $t$, $\hat{u}=(\hat{u}_1,\cdots,\hat{u}_N)$ is an $\epsilon_1$-optimal feedback response,
	$$\frac{1}{N}J^{(N)}_{\rm soc}(\hat{u}_0,\hat{u})\leq \frac{1}{N}J^{(N)}_{\rm soc}(\hat{u}_0,u)+\epsilon_1, \ \hbox{for any}\ u\in \mathcal{U}_c, $$
	where both $\hat{u}_i$ and $u_i$ have the form $Kx_i+\bar{K}\bar{x}+K_0x_0$, $i=1,\ldots N$;
	
	(ii) For any $u_0\in \mathcal{U}_c^0$, $J_0(\hat{u}_0,\hat{u} )\leq J_0(u_0,u)+\epsilon_2$, where $u_0$ has the form $P_0x_0+\bar{P}\bar{x}$ and $\hat{u},u$ are $\epsilon_1$-optimal feedback responses to strategies $\hat{u}_0,u_0$, respectively.

 \end{definition}

From now on, we may suppress the notation of time $t$ if necessary.

%

%
\section{Open-loop Solutions and FBSDEs}
\subsection{The MF Social Control Problem for $N$ Followers}\label{sec3.1}

Suppose $u_0=\check{u}_0\in \mathcal{U}_{d}^0$ 
is fixed. 
We first consider the following  centralized social control problem for $N$ followers.

\textbf{(P1)}: minimize $ {J}_{\rm soc}^{(N)}(u,\check{u}_0)$ over $u\in{\mathcal U}_{c}$, where
{\begin{align*}
  {J}_{\rm soc}^{(N)}(u,\check{u}_0)=&\sum_{i=1}^N
\mathbb{E}\int_0^{T}
\Big\{\big\|x_i
   -\Gamma x^{(N)}-\Gamma_1x_0\big\|^2_{Q}+\|u_i\|^2_{R}+2u_i^TL\check{u}_0+\|\check{u}_0\|^2\Big\}dt.
\end{align*}}

Denote $$\begin{aligned}
&Q_{\Gamma}\stackrel{\Delta}{=}Q\Gamma+\Gamma^TQ-\Gamma^TQ\Gamma,\quad Q_{\Gamma_1}\stackrel{\Delta}{=}Q\Gamma_1-\Gamma^TQ\Gamma_1,\\ &H_{\bar{\Gamma}}\stackrel{\Delta}{=}H\bar{\Gamma}+\bar{\Gamma}^TH-\bar{\Gamma}^TH\bar{\Gamma},\quad  H_{\bar{\Gamma}_1}\stackrel{\Delta}{=}H\bar{\Gamma}_1-\bar{\Gamma}^TH\bar{\Gamma}_1.
\end{aligned}$$

By examining the social cost variation, we obtain the centralized optimal control for followers.
\begin{theorem}\label{thm1}
  Suppose $Q\geq0$ and $R>0$. Then the following system of FBSDEs has a set of solutions $\{x_i,p_i,q_i^j,i,j=0,1,\cdots,N\}$:
  \begin{equation}\label{eq3a}
  \left\{
   \begin{aligned}
   d\check{x}_0=&[A_0\check{x}_0+B_0\check{u}_0
+G_0\check{x}^{(N)}]dt
+D_0dW_0,\cr
  d\check{x}_i=& \big[A\check{x}_i-BR^{-1}B^T\check{p}_i+G\check{x}^{(N)}+F\check{x}_0+(B_1-BR^{-1}L)\check{u}_0 \big]dt+DdW_i,\cr
   d\check{p}_0=&-\big[A^T_0\check{p}_0+F^T\check{p}^{(N)}-Q_{\Gamma_1}^T \check{x}^{(N)}+\Gamma_1^TQ\Gamma_1 \check{x}_0\big]dt
+\sum_{j=0}^N\check{q}_0^jdW_j,\\
  d\check{p}_i=&-\big[A^T\check{p}_i+G^T\check{p}^{(N)}+G^T_0\check{p}_0+Q\check{x}_i-Q_{\Gamma}\check{x}^{(N)}-Q_{\Gamma_1}\check{x}_0\big]dt+\sum_{j=0}^N\check{q}_i^jdW_j,\\
\check{x}_0(0)&=\xi_{0},\ \check{x}_i(0)=\xi_{i},\  \check{p}_0(T)=-H_{\bar{\Gamma}_1} ^T\check{x}^{(N)}(T)+\bar{\Gamma}_1^TH\bar{\Gamma}_1\check{x}_0(T),\ \cr
\check{p}_i(T)&=H\check{x}_i(T)-H_{\bar{\Gamma}}\check{x}^{(N)}(T)-H_{\bar{\Gamma}_1} \check{x}_0(T),\ i=1,\cdots,N.
\end{aligned}\right.
  \end{equation}
Furthermore, the optimal control laws of followers are given by  $\check{u}_i=-R^{-1}(B^T\check{p}_i+ L\check{u}_0), i=1,\cdots,N$.

\end{theorem}

{\it Proof.} See Appendix A. \hfill{$\Box$}

Let $\bar{B}_1\stackrel{\Delta}{=}B_1-BR^{-1}L$. After applying the controls of followers
$\check{u}_i=-R^{-1}(B^T\check{p}_i+L\check{u}_0)$, we have
\begin{equation}
  d\check{x}_i=( A\check{x}_i-BR^{-1}B^T\check{p}_i+G\check{x}^{(N)}+F\check{x}_0+\bar{B}_1\check{u}_0)dt+DdW_i,\ i=1,\ldots,N.
\end{equation}
This with (\ref{eq3a}) leads to
\begin{equation*}
\left\{\begin{aligned}
  d\check{x}^{(N)}=& [(A+G)\check{x}^{(N)}-BR^{-1}B^T\check{p}^{(N)}+F\check{x}_0+\bar{B}_1\check{u}_0]dt+\frac{1}{N}\sum_{i=1}^N\sigma dW_i, \ \check{x}^{(N)}(0)=\xi^{(N)},\cr
   d\check{p}^{(N)}=&-\big[(A+G)^T\check{p}^{(N)}+G^T_0\check{p}_0+(Q-Q_{\Gamma})\check{x}^{(N)}
-(I-\Gamma)^TQ\Gamma_1\check{x}_0\big]dt\\&+\frac{1}{N}\sum_{i=1}^N\sum_{j=0}^N\check{q}_i^jdW_j, \ \check{p}^{(N)}(T)=(H-H_{\bar{\Gamma}})x^{(N)}(T)-H_{\bar{\Gamma}_1} x_0(T),
\end{aligned}\right.
\end{equation*}
where $\xi^{(N)}=\frac{1}{N}\sum_{i=1}^N\xi_i$. Let $N\to\infty$. By the law of large numbers, we may approximate $\check{x}_0,\check{x}^{(N)}$, $\check{p}_0,\check{p}^{(N)}$ by $\bar{x}_0,\bar{x}$, $\bar{p}_0,\bar{p}$, which satisfy
\begin{equation}\label{eq4a}
\left\{\begin{aligned}
 d\bar{x}_0=&(A_0\bar{x}_0+B_0\check{u}_0
+G_0\bar{x})dt
+D_0dW_0,\ \bar{x}_0(0)=\xi_0\cr
  d\bar{x}=& [(A+G)\bar{x}-BR^{-1}B^T\bar{p}+F\bar{x}_0+\bar{B}_1\check{u}_0]dt,\ \bar{x}(0)=\bar{\xi},\cr
     d\bar{p}_0=&-\big(A^T_0\bar{p}_0+F^T\bar{p}-Q_{{\Gamma}_1} ^T\bar{x} +\Gamma_1^TQ\Gamma_1 \bar{x}_0\big)dt
+\bar{q}_0^0dW_0,\ \bar{p}_0(T)=-H_{\bar{\Gamma}_1} ^T\bar{x}(T)+\bar{\Gamma}_1^TH\bar{\Gamma}_1\bar{x}_0(T),\\
   d\bar{p}=&-\big[(A+G)^T\bar{p}+G^T_0\bar{p}_0+(Q-Q_{\Gamma})\bar{x}
-Q_{{\Gamma}_1}\bar{x}_0\big]dt\\&+\bar{q}^0dW_0,\
\bar{p}(T)=(H-H_{\bar{\Gamma}})\bar{x}(T)-H_{\bar{\Gamma}_1} \bar{x}_0(T).
\end{aligned}\right.
\end{equation}
Based on this with (\ref{eq3a}), we construct the following FBSDEs
\begin{equation}\label{eq11-a}
\left\{\begin{aligned}
  d\bar{x}_i=& \big(A\bar{x}_i+G\bar{x}-BR^{-1}B^T\bar{p}_i+F\bar{x}_0+\bar{B}_1\check{u}_0\big)dt+DdW_i,\ \bar{x}_i(0)=\xi_i,\cr
        d\bar{p}_i=&-\big(A^T\bar{p}_i+G^T\bar{p}+G^T_0\bar{p}_0+Q\bar{x}_i-Q_{\Gamma}\bar{x}
-Q_{\Gamma_1}\bar{x}_0\big)dt+\bar{q}_i^idW_i\\
&+\bar{q}_i^0dW_0,\ \bar{p}_i(T)=H\bar{x}_i(T)-H_{\bar{\Gamma}}\bar{x}(T)-H_{\bar{\Gamma}_1} \bar{x}_0(T),
\end{aligned}\right.
\end{equation}
and the decentralized control laws of followers are given by
\begin{equation}\label{eq13a}
  u_i^*=-R^{-1}(B^T\bar{p}_i+L\check{u}_0),\ i=1,\cdots,N.
\end{equation}
\begin{remark}
  By (\ref{eq4a})-(\ref{eq11-a}), we have $\bar{x}=\mathbb{E}[\bar{x}_i|\mathcal{F}^0]$ and $\bar{p}=\mathbb{E}[\bar{p}_i|\mathcal{F}^0]$. Thus, (\ref{eq11-a}) is equivalent to
  \begin{equation*}
\left\{\begin{aligned}
  d\bar{x}_i=& \big[A\bar{x}_i+G\mathbb{E}[\bar{x}_i|\mathcal{F}^0]-BR^{-1}B^T\bar{p}_i+F\bar{x}_0+\bar{B}_1\check{u}_0\big]dt+DdW_i,\ \bar{x}_i(0)=\xi_i,\cr
        d\bar{p}_i=&-\big[A^T\bar{p}_i+G^T\mathbb{E}[\bar{p}_i|\mathcal{F}^0]+G^T_0\bar{p}_0+Q\bar{x}_i
        -Q_{\Gamma}\mathbb{E}[\bar{x}_i|\mathcal{F}^0]
-Q_{\Gamma_1}\bar{x}_0\big]dt\\&+\bar{q}_i^idW_i+\bar{q}_i^0dW_0,\ \bar{p}_i(T)=H\bar{x}_i(T)-H_{\bar{\Gamma}}\mathbb{E}[\bar{x}_i(T)|\mathcal{F}^0]-H_{\bar{\Gamma}_1} \bar{x}_0(T),
\end{aligned}\right.
\end{equation*}
which is a (conditional) MF FBSDE \cite{Y13}.
\end{remark}

We now use the idea inspired by \cite{YZ99}, \cite{ZQ16}, \cite{WZZ20} to decouple the FBSDEs
 (\ref{eq4a}) and (\ref{eq11-a}).
Let $\bar{p}_0=\Pi_0 \bar{x}+M_0\bar{x}_0+\bar{\varphi}_0,\ \bar{p}=\bar{\Pi} \bar{x}+M\bar{x}_0+\bar{\varphi}.$
 Then applying It\^{o}'s formula, we obtain
 $$\begin{aligned}
  d\bar{p}_0=&(\dot{{\Pi}}_0\bar{x}+\dot{M}_0\bar{x}_0)dt+d\bar{\varphi}_0+{\Pi}_0\big[(A+G){\bar{x}}
  -B{R^{-1}}B^T(\bar{\Pi}\bar{x}+M\bar{x}_0+\bar{\varphi})
+F\bar{x}_0\cr
&+(B_1-BR^{-1}L)\check{u}_0\big]dt+M_0\big[(A_0\bar{x}_0+B_0\check{u}_0+G_0\bar{x})dt+D_0dW_0\big]\cr
=&-\big[A^T_0(\Pi_0 \bar{x}+M_0\bar{x}_0+\bar{\varphi}_0)+F^T(\bar{\Pi}\bar{x}+M\bar{x}_0+\bar{\varphi})\cr
&-Q_{\Gamma_1}^T\bar{x}+\Gamma_1^TQ\Gamma_1 \bar{x}_0\big]
+\bar{q}_0^0dW_0,
\end{aligned}
$$
which implies
\begin{align}
 \label{eq6a}
&\dot{\Pi}_0+\Pi_0(A+G)+A_0^T\Pi_0-\Pi_0 BR^{-1}B^T\bar{\Pi}+M_0G_0+F^T\bar{\Pi}-Q_{\Gamma_1}^T=0,\  {\Pi}_0(T)=
-H_{\bar{\Gamma}_1}^T,
\\
 \label{eq6ab}&\dot{M}_0+M_0A_0+A_0^TM_0+(F^T-\Pi_0BR^{-1}B^T)M+\Pi_0F+\Gamma_1^TQ\Gamma_1=0,\ M_0(T)=\bar{\Gamma}_1^TH\bar{\Gamma}_1,\\
 \label{eq7a}
&d\bar{\varphi}_0=-\big[A_0^T\bar{\varphi}_0+(F^T-\Pi_0BR^{-1}B^T)\bar{\varphi}+(\Pi_0\bar{B}_1+M_0B_0)\check{u}_0\big]dt\cr
&+(\bar{q}_0^0-M_0D_0)dW_0,\  \bar{\varphi}_0(T)=0.
\end{align}
Besides, by It\^{o}'s formula, we have
$$\begin{aligned}
d\bar{p}=\ &(\dot{\bar{\Pi}}\bar{x}+\dot{M}\bar{x}_0)dt+\bar{\Pi}\big[(A+G)\bar{x}-B{R^{-1}}B^T(\bar{\Pi}\bar{x}+M\bar{x}_0+\bar{\varphi})+F\bar{x}_0\cr&
+\bar{B}_1\check{u}_0\big]dt
+M\big[(A_0\bar{x}_0+B_0\check{u}_0+G_0\bar{x})dt+D_0dW_0\big]dt+d\bar{\varphi}\cr
=&-\big[(A+G)^T(\bar{\Pi}\bar{x}+M\bar{x}_0+\bar{\varphi})+G^T_0({\Pi}_0\bar{x}+M_0\bar{x}_0+\bar{\varphi}_0)+(Q-Q_{\Gamma})\bar{x}\\&
-Q_{\Gamma_1}\bar{x}_0\big]dt+\bar{q}^0dW_0,
\end{aligned}$$
which implies
\begin{align}\label{eq6}
&\dot{\bar{\Pi}}+(A+G)^T\bar{\Pi}+\bar{\Pi}(A+G)-\bar{\Pi} BR^{-1}B^T\bar{\Pi}\cr&
+MG_0 +G_0^T\Pi_0+Q-Q_{\Gamma}=0,\ \bar{\Pi}(T)=H-H_{\bar{\Gamma}},\\
\label{eq6-b}
&\dot{M}+(A+G)^TM+MA_0-\bar{\Pi}BR^{-1}B^TM+G_0^TM_0\cr
&+\bar{\Pi}F-Q_{\Gamma_1}=0, \ M(T)=-H_{\bar{\Gamma}_1}, \\
 \label{eq7}
&d\bar{\varphi}=-\big[(A+G-BR^{-1}B^T\bar{\Pi})^T\bar{\varphi}+G^T_0\bar{\varphi}_0
+(\bar{\Pi}\bar{B}_1+MB_0)\check{u}_0\big]dt\cr&+(\bar{q}^0-MD_0)dW_0,\  \bar{\varphi}(T)=0.
\end{align}
By observing (\ref{eq6a}), (\ref{eq6ab}), (\ref{eq6}) and (\ref{eq6-b}), we have $\Pi_0=M^T$, and $M_0,\bar{\Pi}$ are symmetric.

From (\ref{eq4a}) and (\ref{eq11-a}),
\begin{equation}\label{eq17-a}
\left\{\begin{aligned}
  d(\bar{x}_i-\bar{x})=&\big[A(\bar{x}_i-\bar{x})-BR^{-1}B(\bar{p}_i-\bar{p})\big]dt+DdW_i,\ \bar{x}_i(0)-\bar{x}(0)=\xi_i-\bar{\xi},\cr
        d(\bar{p}_i-\bar{p})=&-\big[A^T(\bar{p}_i-\bar{p})+Q(\bar{x}_i-\bar{x})
\big]dt+\bar{q}_i^idW_i+(\bar{q}_i^0-\bar{q}^0)dW_0,\ \cr
\bar{p}_i(T)-\bar{p}(&T)=H(\bar{x}_i(T)-\bar{x}(T)).
\end{aligned}\right.
\end{equation}
Suppose $\bar{p}_i-\bar{p}=\Pi(\bar{x}_i-\bar{x})$. By It\^{o}'s formula,
$$d(\bar{p}_i-\bar{p})=\dot{\Pi}(\bar{x}_i-\bar{x})dt+\Pi\big[\big(A(\bar{x}_i-\bar{x})-BR^{-1}B(\bar{p}_i-\bar{p})\big)dt+DdW_i\big].$$
Comparing this with (\ref{eq17-a}), it follows that
 $\bar{q}_i^i=\Pi D$, $\bar{q}_i^0=\bar{q}^0$, and $\Pi$ should satisfy
\begin{equation}
  \dot{{\Pi}}+ A^T{\Pi}+{\Pi}A-{\Pi}BR^{-1}B^T{\Pi}+Q=0,\ {\Pi}(T)=H. \label{eq5}
\end{equation}

\begin{theorem}
  If the equations (\ref{eq6a}), (\ref{eq6}) and (\ref{eq5}) admit a set of solutions, then the decentralized control law (\ref{eq13a}) has a feedback representation:
  \begin{align}\label{eq-d1}
    {u}^*_i=-R^{-1}B^T[\Pi\bar{x}_i+(\bar{\Pi}-\Pi)\bar{x}+M\bar{x}_0+\bar{\varphi}]-R^{-1}L\check{u}_0,
  \end{align}
  where  $\Pi,\bar{\Pi},M,\bar{\varphi},\bar{x}_i, \bar{x},\bar{x}_0$ are determined by (\ref{eq5}), (\ref{eq6})-
(\ref{eq7}), (\ref{eq4a}) and (\ref{eq11-a}).
\end{theorem}
\emph{Proof.} Note that (\ref{eq7a}) and (\ref{eq7}) are linear backward SDEs, and hence both admit a solution, respectively. By the above discussion,
$\bar{p}_i=\bar{p}+\bar{p}_i-\bar{p}
=\Pi\bar{x}_i+(\bar{\Pi}-\Pi)\bar{x}+M\bar{x}_0+\bar{\varphi}.
$
\hfill{$\Box$}
\begin{remark}
  Note that $Q-Q_{\bar{\Gamma}}  \geq0$ and  $H-H_{\bar{\Gamma}}
  \geq0$. If $G_0=0$, then by \cite{YZ99}, Riccati equation (\ref{eq6}) admits a unique solution $\bar{\Pi}\geq 0$. Accordingly, (\ref{eq6a})-(\ref{eq7a}) and (\ref{eq6-b})-(\ref{eq7}) have solutions since they are linear matrix (vector) differential equations.
\end{remark}


\subsection{The Optimization Problem for the Leader}

After applying the controls of followers 
 in (\ref{eq-d1}), we have an optimal control problem for the leader.

\textbf{(P2)}: minimize $ {J}_{0}(u_0,{u}^*)$ over $u_0\in L_{\mathcal{F}_t}^2(0,T;\mathbb{R}^m)$, where
\begin{align}
{J}_0(u_0,{u}^*)=&\mathbb{E}\int_{0}^{T}\big[\|x_0- \Gamma_0{x}_*^{(N)}\|_{Q_0}^{2}+\|u_0\|_{R_0}^{2}\big]dt+\mathbb{E}\big[\|x_0(T)-\bar{\Gamma}_0x_*^{(N)}(T)\|^2_{H_0}\big],\\ dx_0=&[A_0x_0+B_0u_0
+G_0{x}_*^{(N)}]dt
+D_0dW_0,\ x_0(0)=\xi_0,\label{eq18d}\\
d{x}_i^*=&\big[A{x}_i^*+G{x}_*^{(N)}-BR^{-1}B^T\big(\Pi \bar{x}_i+(\bar{\Pi}-\Pi) \bar{x}+M\bar{x}_0+\bar{\varphi}\big)\cr
&+Fx_0+\bar{B}_1 u_0\big]dt+DdW_i,\ {x}^*_i(0)=\xi_i,\label{eq19d}
\end{align}
where $\bar{x}_i, \bar{x},\bar{\varphi}$ satisfy (\ref{eq4a}), (\ref{eq11-a}) and  (\ref{eq7}), respectively.
 In the above, $x_i^*$ is the realized state under the control $u^*_i, i=1,\cdots,N$ and ${x}_*^{(N)}=\frac{1}{N}\sum_{i=1}^N x_i^*$. From (\ref{eq19d}), we have
 \begin{align}\label{eq23}
   d{x}_*^{(N)}=&\big[(A+G){x}_*^{(N)}-BR^{-1}B^T\big(\Pi \bar{x}^{(N)}+(\bar{\Pi}-\Pi) \bar{x}+M\bar{x}_0+\bar{\varphi}\big)\cr
&+Fx_0+\bar{B}_1 u_0\big]dt+\frac{1}{N}\sum_{i=1}^NDdW_i,\ {x}^*_i(0)=\xi_i,
 \end{align}
where $\bar{x}^{(N)}=\frac{1}{N}\sum_{i=1}^N\bar{x}_i$.
 Since $\{W_i\}$ are independent Wiener processes and $\{x_i(0)\}$ are independent random variables, for the large $N$, it is plausible to replace $\bar{x}^{(N)},{x}_*^{(N)}$ by $\bar{x}$, which evolves by (\ref{eq7}). 
Then we have the following optimal control problem for the leader.

\textbf{(P2$^\prime$)}: minimize $ \bar{J}_{0}(u_0)$ over $u_0\in L_{\mathcal{F}^0}^2(0,T;\mathbb{R}^m)$, where
\begin{align}
\bar{J}_0(u_0)=&\mathbb{E}\int_{0}^{T}\big[\|\bar{x}_0- \Gamma_0\bar{x}\|_{Q_0}^{2}+\|u_0\|_{R_0}^{2}\big]dt+\mathbb{E}\big[\|\bar{x}_0(T)-\bar{\Gamma}_0\bar{x}(T)\|^2_{H_0}\big],\cr
d\bar{x}_0=&[A_0\bar{x}_0+B_0u_0+G_0\bar{x}]dt
+D_0dW_0,\ \bar{x}_0(0)=\xi_0,\label{eq8b}
\\ \label{eq8c}
d\bar{x}=&[(A+G-BR^{-1}B^T\bar{\Pi}) \bar{x}+(F-BR^{-1}B^TM)\bar{x}_0+\bar{B}_1 u_0-BR^{-1}B^T\bar{\varphi}]dt,\ \bar{x}(0)=\bar{\xi},\\  \label{eq8d}
d\bar{\varphi}_0=&-\big[A_0^T\bar{\varphi}_0+(F^T-\Pi_0 BR^{-1}B^T)\bar{\varphi}+(\Pi_0\bar{B}_1+M_0B_0)u_0\big]dt\cr
&+(\bar{q}_0^0-M_0D_0)dW_0,\ \bar{\varphi}_0(T)=0,\\ \label{eq8e}
d\bar{\varphi}=&-\big\{(A+G-BR^{-1}B^T\bar{\Pi})^T\bar{\varphi}+G^T_0\bar{\varphi}_0
\cr&+(\bar{\Pi}\bar{B}_1+MB_0)u_0\big\}dt+(\bar{q}^0-MD_0)dW_0,\  \bar{\varphi}(T)=0.
\end{align}


Define the FBSDE
  \begin{equation}\label{eq16}
  \left\{
\begin{aligned}
d\bar{x}_0=&\big\{A_0\bar{x}_0-B_0R^{-1}_0[B_0^Ty_0+\bar{B}_1^T\bar{y}+({\Pi}_0\bar{B}_1+M_0B_0)^T\psi_0+(\bar{\Pi}\bar{B}_1+MB_0)^T\psi]\cr
&+G_0\bar{x}\big\}dt
+D_0dW_0,\ \bar{x}_0(0)=\xi_0,\\
d\bar{x}=&[(A+G-BR^{-1}B^T\bar{\Pi}) \bar{x}+(F-BR^{-1}B^TM)\bar{x}_0-BR^{-1}B^T\bar{\varphi}\cr
&-\bar{B}_1R^{-1}_0[B_0^Ty_0+\bar{B}_1^T\bar{y}+({\Pi}_0\bar{B}_1+M_0B_0)^T\psi_0+(\bar{\Pi}\bar{B}_1+MB_0)^T\psi]dt,\ \bar{x}(0)=\bar{\xi},\\
d{\bar{\varphi}_0}=&-\big\{A_0^T\varphi_0+(F-BR^{-1}B^TM)^T\varphi-(\Pi_0\bar{B}_1+M_0B_0)R^{-1}_0[B_0^Ty_0+\bar{B}_1^T\bar{y}\cr
&+({\Pi}_0\bar{B}_1+M_0B_0)^T\psi_0+(\bar{\Pi}\bar{B}_1+MB_0)^T\psi]
\big\}dt+(\bar{q}_0^0-M_0D_0)dW_0,\  \bar{\varphi}_0(T)=0,\cr
d{\bar{\varphi}}=&-\big\{(A+G-BR^{-1}B^T\bar{\Pi})^T\varphi+G^T_0\varphi_0
-(\bar{\Pi}\bar{B}_1+MB_0)R^{-1}_0[B_0^Ty_0+\bar{B}_1^T\bar{y}\cr
&+({\Pi}_0\bar{B}_1+M_0B_0)^T\psi_0+(\bar{\Pi}\bar{B}_1+MB_0)^T\psi]
\big\}dt+(\bar{q}^0-MD_0)dW_0,\  \varphi(T)=0,\cr
  dy_0=&-\big[ A^T y_0+(F-BR^{-1}B^TM)^T\bar{y}
  +Q_0(\bar{x}_0-\Gamma_0\bar{x})\big]dt+\beta_0dW_0,\ \cr
&\qquad y_0(T)=H_0(\bar{x}_0(T)-\bar{\Gamma}_0\bar{x}(T)),\\
  d\bar{y}=& -\big[(A+G-BR^{-1}B^T\bar{\Pi})^T \bar{y}+G_0^Ty_0-\Gamma_0^TQ_0 (\bar{x}_0-\Gamma_0\bar{x})\big]dt+\bar{\beta} dW_0,\cr
   &\qquad \bar{y}(T)=-\bar{\Gamma}_0^TH_0(\bar{x}_0(T)-\bar{\Gamma}_0\bar{x}(T)),\\
  d\psi_0=&\ (A_0\psi_0+G_0\psi) dt,\ \psi_0(0)=0,\\
   d\psi=&\ [(A+G-BR^{-1}B^T\bar{\Pi})\psi+(F-BR^{-1}B^T\Pi_0)\psi_0-BR^{-1}B^T \bar{y}]dt,\ \psi(0)=0.
\end{aligned}
\right.
\end{equation}

\begin{theorem}\label{thm3.3}
  Assume A1)-A2) hold, and (\ref{eq16}) has a solution over $[0,T]$. Then Problem (P2$^{\prime}$) admits an optimal control $$u_0^*=
  -R^{-1}_0[B_0^Ty_0+\bar{B}_1^T\bar{y}+({\Pi}_0\bar{B}_1+M_0B_0)^T\psi_0+(\bar{\Pi}\bar{B}_1+MB_0)^T\psi],$$ where
  $y_0$, $\bar{y}$, $\psi_0$ and ${\psi}$ satisfy (\ref{eq16}).
\end{theorem}

\emph{Proof.} See Appendix A. \hfill{$\Box$}

Denote $X=[\bar{x}_0^T,\bar{x}^T,\psi_0^T,\psi^T
]^T,
Y=[y_0^T,\bar{y}^T,\bar{\varphi}_0^T,\bar{\varphi}^T]^T$, $Z=[\beta_0^T,\bar{\beta}^T,(\bar{q}_0^0)^T,(\bar{q}_i^0)^T]^T$. Let
$$
\mathcal{A}=\left[\begin{array}{cccc}
A_0&G_0&-B_0R_0^{-1}\Xi_0^T&-B_0R_0^{-1}\bar{\Xi}^T\\
F-BR^{-1}B^TM&\hat{A}&-\bar{B}_1R_0^{-1}\Xi_0^T&-\bar{B}_1R_0^{-1}\bar{\Xi}^T\\
0&0&A_0&G_0\\
0&0&F-BR^{-1}B^T\Pi_0&\hat{A}
\end{array}\right],
$$
$$\mathcal{B}=\left[\begin{array}{cccc}
B_0R_0^{-1}B_0^T&B_0R_0^{-1}\bar{B}_1^T&0&0\\
\bar{B}_1R_0^{-1}B_0^T&\bar{B}_1R_0^{-1}\bar{B}_1^T&0&BR^{-1}B^T\\
0&0&0&0\\
0&BR^{-1}B^T&0&0
\end{array}\right],\
\mathcal{H}_0=\left[\begin{array}{cccc}
H_0&-H_0\bar{\Gamma}_0&0&0\cr
-\bar{\Gamma}_0^TH_0&\bar{\Gamma}_0^TH_0\bar{\Gamma}_0&0&0\cr
0&0&0&0\cr
0&0&0&0
\end{array}\right],$$
$$\mathcal{Q}=\left[\begin{array}{cccc}
-Q_0&Q_0\Gamma_0&0&0\\
\Gamma_0^TQ_0&-\Gamma_0^TQ_0\Gamma_0&0&0\\
0&0&\Xi_0R_0^{-1}\Xi_0^T&\Xi_0R_0^{-1}\bar{\Xi}^T\\
0&0&\bar{\Xi}R_0^{-1}\Xi_0^T&\bar{\Xi}R_0^{-1}\bar{\Xi}^T
\end{array}\right],$$
$\mathcal{D}_0=[D_0^T,0,0,0]^T$  and  $\bar{\mathcal{D}}_0=[0,0,D_0^TM_0^T,D_0^TM^T]^T$, where
 $\hat{A}\stackrel{\Delta}{=}A+G-BR^{-1}B^T\bar{\Pi}$, $\Xi_0\stackrel{\Delta}{=}{\Pi}_0\bar{B}_1+M_0B_0$ and $\bar{\Xi}\stackrel{\Delta}{=}\bar{\Pi}\bar{B}_1+MB_0$.
With the above notions, we can rewrite (\ref{eq16}) as
\begin{equation}
\label{eq17}
  \left\{
\begin{aligned}
  dX&=(\mathcal{A}X-\mathcal{B}Y)dt+\mathcal{D}_0dW_0,\ X(0)=[\xi_0^T, \bar{\xi}^T,0,0]^T,\\
  dY&=(\mathcal{Q}X-\mathcal{A}^TY)dt+(Z-\bar{\mathcal{D}}_0)dW_0,\ Y(T)=\mathcal{H}_0X(T).
\end{aligned}
\right.
\end{equation}

For further analysis,  assume

\textbf{(A3)} FBSDE (\ref{eq17}) admits a solution $(X,Y,Z)$ over $[0,T]$.

We now provide a sufficient condition to guarantee (A3).
\begin{proposition}
 If  the equation
 \begin{align}\label{eq18}
 &\dot{\mathcal{P}}+\mathcal{P}\mathcal{A}+\mathcal{A}^T\mathcal{P}-\mathcal{P}\mathcal{B}\mathcal{P}-\mathcal{Q}=0,
\quad \mathcal{P}(T)=\mathcal{H}_0
  \end{align}
  has a solution in $[0,T]$, then (A3) holds.
\end{proposition}
\emph{Proof.}
Let $Y=\mathcal{P}X$. Then, we have
\begin{align*}
  dY&=\dot{\mathcal{P}}Xdt
  +\mathcal{P}[(\mathcal{A}X-\mathcal{B}\mathcal{P}X)dt+\mathcal{D}_0d{W}_0]\\
&=(\mathcal{Q}X-\mathcal{A}^T\mathcal{P}X)dt+ZdW_0.
\end{align*}
This implies that $Z=\mathcal{PD}_0+\bar{\mathcal{D}}_0$ and $\mathcal{P}$ should satisfy (\ref{eq18}).
  If (\ref{eq18}) has a solution in $[0,T]$, then by \cite{my1999}, FBSDE (\ref{eq16}) admits an adapted solution. \hfill$\Box$
\begin{remark}
 Note that $\mathcal{B}$, $\mathcal{Q}$  and $\mathcal{H}_0$ are symmetric matrices. We find that  (\ref{eq18}) is a symmetric Riccati equation.
 The existence condition of its solution may be referred in \cite{AF03} and \cite{my1999}.
\end{remark}

Under (A3), we design the following decentralized control laws
 \begin{equation}\label{eq27c}
 \left\{
   \begin{aligned}
     {u}_0^*=&-R_0^{-1}[B_0^Ty_0+\bar{B}_1^T\bar{y}+({\Pi}_0\bar{B}_1+M_0B_0)^T\psi_0+(\bar{\Pi}\bar{B}_1+MB_0)^T\psi],\\
     {u}_i^*=&-R^{-1}B^T[\Pi\bar{x}_i+(\bar{\Pi}-\Pi)\bar{x}+M\bar{x}_0+\bar{\varphi}]-R^{-1}Lu_0^*,
   \end{aligned}
   \right.
 \end{equation}
 where $y_0,\bar{y},\psi_0,\psi, \bar{x}_0,\bar{x},\bar{\varphi}$ are given by (\ref{eq16}), and $\bar{x}_i$ satisfies
\begin{align}\label{eq40}
d\bar{x}_i=& \big[(A-BR^{-1}B^T\Pi)\bar{x}_i+(G-BR^{-1}B^T(\bar{\Pi}-\Pi))\bar{x}-BR^{-1}B^T\bar{\varphi}\cr
&+(F-BR^{-1}B^TM)\bar{x}_0\big]dt+DdW_i,\ \bar{x}_i(0)=\xi_i.
\end{align}
%
%
%

 \begin{theorem}\label{thm3.4}
   For Problem (\ref{eq1})-(\ref{eq4}), assume that (A1)-(A3) hold. Then $({u}_0^*,{u}_1^*,\cdots,
{u}_N^*)$ given in (\ref{eq27c}) is an $(\varepsilon_1,\varepsilon_2)$-Stackelberg equilibrium, 
where $\varepsilon_i=O(1/\sqrt{N})$, $i=1,2$.
 \end{theorem}

 \emph{Proof.} See Appendix \ref{app2}. \hfill{$\Box$}

\begin{theorem}\label{thm3.5}
   Assume that (A1)-(A3) hold, and $\xi_i,i=1,\cdots,N$ have the same variance. Then under the decentralized control (\ref{eq27c}), the asymptotic average social cost of followers is given by
  \begin{equation}
    \label{eq46a}
\lim_{N\to\infty}\frac{1}{N}J_{\rm soc}^{(N)}(u^*,u_0^*)=\mathbb{E}\big[\|\xi_i\|^2_{\Pi(0)}+\|\bar{\xi}\|^2_{\bar{\Pi}(0)
-\Pi(0)}+\|\xi_0\|^2_{M_0(0)}+2\bar{\xi}^T\Pi_0(0)\xi_0\big]+m_T,
  \end{equation}
    and the asymptotic cost of the leader is given by
 \begin{equation}\label{eq46}
   \lim_{N\to\infty}J_{0}(u^*,u_0^*)=\mathbb{E}\big[{\xi}^T_0y_0(0)+\bar{\xi}^T\bar{y}(0)+\int_0^T(D_0^T\beta_0)dt\big],
 \end{equation}
  where
\begin{align*}
  m_T\stackrel{\Delta}{=}& \mathbb{E}\big[2\bar{\xi}^T\bar{\varphi}(0)
+2{\xi}^T_0\bar{\varphi}_0(0)  \big]+\mathbb{E}\int_0^T\big[2\bar{\varphi}_0^TB_0{u}_0^*-2\bar{\varphi}^TBR^{-1}L{u}_0^*-\|B^T\bar{\varphi}\|^2_{R^{-1}}\cr
&+\|{u}_0^*\|^2_{R_{1}}-\|L{u}_0^*\|^2_{R^{-1}}+D^T\Pi D+D_0^TM_0D_0+D_0^T(\bar{q}_0^0-M_0D_0)\big]dt.
\end{align*}

\end{theorem}

 \emph{Proof.} See Appendix \ref{app3}. \hfill{$\Box$}

\section{Feedback Solutions and Riccati Equations}

In this section, we consider the feedback solutions to the leader-follower MF  game.

\subsection{MF  Social Control for $N$ Followers}

Similar to 
the line of Section \ref{sec3.1},
we first consider a centralized social control problem for $N$ followers.
Thus, the leader is presumed to adopt a feedback control with the following form
\begin{equation}
  \label{eq35}
u_0=P_0x_0+\bar{P}x^{(N)},
\end{equation}
 where $P_0$ and $\bar{P}$ are fixed. 
This leads to a centralized social control problem for $N$ followers.

\textbf{(P3)}: minimize $ {J}_{\rm soc}^{(N)}(u)$ over $u\in{\mathcal U}_{c}$, where $u_0=P_0x_0+\bar{P}x^{(N)}$ and
{\begin{align*}
  {J}_{\rm soc}^{(N)}(u)=&\sum_{i=1}^N
\mathbb{E}\int_0^{T}
\Big\{\big\|x_i
   -\Gamma x^{(N)}-\Gamma_1x_0\big\|^2_{Q}+\|u_i\|^2_{R}+2u_i^TL(P_0x_0+\bar{P}x^{(N)})\cr
   &+\|P_0x_0+\bar{P}x^{(N)}\|^2_{R_1}\Big\}dt+\mathbb{E}
\big[\big\|x_i(T)
   -\bar{\Gamma} x^{(N)}(T)-\bar{\Gamma}_1x_0(T)\big\|^2_{H}\big].
\end{align*}}

An additional assumption is now introduced.

\textbf{(A4)} $R-LR_1^{-1}L^T\geq 0$.

By examining the social cost variation, we obtain the optimal control laws for $N$ followers.
\begin{theorem}\label{thm4.1}
  Suppose that (A2) and (A4) hold. 
  Then the following system of FBSDEs admits a set of adapted solutions $\{x_i,p_i,q_i^j,i,j=0,1,\cdots,N\}$:
  \begin{equation}\label{eq3d}
  \left\{
   \begin{aligned}
   dx_0=&\big[A_0x_0+B_0(P_0x_0+\bar{P}x^{(N)})
+G_0x^{(N)}\big]dt
+D_0dW_0,\cr
  dx_i=& [Ax_i+
  B\breve{u}_i+Gx^{(N)}+Fx_0+{B}_1(P_0x_0+\bar{P}x^{(N)})]dt+DdW_i,\cr
   dp_0=&-\big[(A_0+B_0P_0)^Tp_0+(F+B_1P_0)^Tp^{(N)}-Q_{\Gamma_1}^T x^{(N)}+\Gamma_1^TQ\Gamma_1 x_0\cr
   &+P_0^TR_1(P_0x_0+\bar{P}x^{(N)})+P_0^TL^T\breve{u}^{(N)}\big]
+\sum_{j=0}^Nq_0^jdW_j,\\
  dp_i=&-\big[A^Tp_i+(G+B_1\bar{P})^Tp^{(N)}+(G_0+B_0\bar{P})^Tp_0+Qx_i+(\bar{P}^TR_1\bar{P}-Q_{\Gamma})x^{(N)}\\&
+\big(\bar{P}^TR_1P_0-Q_{\Gamma_1}\big)x_0+\bar{P}^TL^T\breve{u}^{(N)}\big]dt+\sum_{j=0}^Nq_i^jdW_j,\\
x_0(0)&=\xi_{0},\ x_i(0)=\xi_{i},\   {p}_0(T)=-H_{\bar{\Gamma}_1} ^Tx^{(N)}(T)+\bar{\Gamma}_1^TH\bar{\Gamma}_1x_0(T),\ \cr
{p}_i(T)&=Hx_i(T)-H_{\bar{\Gamma}}x^{(N)}(T)-H_{\bar{\Gamma}_1} x_0(T),\ i=1,\cdots,N.
\end{aligned}\right.
  \end{equation}
where $\breve{u}^{(N)}=\frac{1}{N}\sum_{i=1}^N\breve{u}_i$ and the optimal strategies of followers are given by
$$\breve{u}_i=-R^{-1}B^Tp_i-R^{-1}L(P_0x_0+\bar{P}x^{(N)}),\quad  i=1,\cdots,N.$$

\end{theorem}

 \emph{Proof.} See Appendix \ref{app4}. \hfill{$\Box$}

By (\ref{eq3d}), we have
\begin{align}\label{eq63}
 d{x}^{(N)}=& \big[(A+G+\bar{B}_1\bar{P}){x}^{(N)}-B{R^{-1}}B^T((K+\bar{K})x^{(N)}+K_0x_0)\cr
&+(F+\bar{B}_1{P}_0)x_0\big]dt+\frac{1}{N}\sum_{i=1}^ND dW_i.
\end{align}
Let $$p_0=\Lambda_0x_0+\bar{\Lambda}x^{(N)},\  p_i=K x_i+\bar{K}x^{(N)}+K_0x_0,\ i=1,\cdots,N.$$ Then applying It\^{o}'s formula, we obtain
$$\begin{aligned}
dp_i=\ &\dot{K}{x}_idt+K\Big[\big(A{x}_i-B{R^{-1}}B^T(Kx_i+\bar{K}x^{(N)}+K_0x_0)+(G+\bar{B}_1\bar{P})x^{(N)}\\
&+(F+\bar{B}_1{P}_0)x_0\big)dt+DdW_i\Big]+\dot{\bar{K}}x^{(N)}dt+\bar{K}\Big\{\big[(A+G+\bar{B}_1\bar{P}){x}^{(N)}\\
&-B{R^{-1}}B^T((K+\bar{K})x^{(N)}+K_0x_0)
+(F+\bar{B}_1{P}_0)x_0\big]dt+\frac{1}{N}\sum_{i=1}^ND dW_i\Big\}\\
&+\dot{K}_0x_0+K_0\big\{[(A_0+B_0P_0)x_0+(G_0+B_0\bar{P})x^{(N)}]dt+D_0dW_0\big\}\cr
\end{aligned}$$
By comparing this with (\ref{eq3d}), we obtain $q_i^0={K}_0D_0$, $q_i^i=\frac{1}{N}\bar{K}D+KD$, $q_i^j=\frac{1}{N}\bar{K}D, \ j\not=i, i,j=1,\cdots,N$, and
\begin{align}\label{eq5b}
 &\dot{K}+ A^TK+KA-KBR^{-1}B^TK+Q=0,\ K(T)=H,\\ \label{eq6b}
&\dot{\bar{K}}+(A+G+\bar{B}_1\bar{P})^T\bar{K}+\bar{K}(A+G+\bar{B}_1\bar{P})-KBR^{-1}B^T\bar{K}-\bar{K}BR^{-1}B^TK\cr
&-\bar{K}BR^{-1}B^T\bar{K}+(G+\bar{B}_1\bar{P})^TK+K(G+\bar{B}_1\bar{P})+K_0(G_0+B_0\bar{P})
\cr
&+(G_0+B_0\bar{P})^T\bar{\Lambda}+\bar{P}^TR_1\bar{P}-Q_{\Gamma}-\bar{P}^TL^TR^{-1}L\bar{P}=0,\ \bar{K}(T)=-H_{\bar{\Gamma}},\\
&\dot{K}_0+(A+G+\bar{B}_1\bar{P})^TK_0+K_0(A_0+B_0P_0)-(K+\bar{K})BR^{-1}B^TK_0+(G_0+B_0\bar{P})^T\Lambda_0\cr&
+(K+\bar{K})(F+\bar{B}_1P_0)- Q_{\Gamma_1}+\bar{P}^TR_1P_0-\bar{P}^TL^TR^{-1}LP_0=0,  {K}_0(T)=\!-H_{\bar{\Gamma}_1}\!.
\end{align}
By It\^{o}'s formula, for any $i=1,\cdots,N$, we have
$$\begin{aligned}
dp_0=\ &\dot{\Lambda}_0{x}_0dt+\Lambda_0\big\{[(A_0+B_0P_0)x_0+(G_0+B_0\bar{P})x^{(N)}]dt+D_0dW_0\big\}+\dot{\bar{\Lambda}}x^{(N)}dt\\&
+\bar{\Lambda}\Big\{\big[(A+G+\bar{B}_1\bar{P}){x}^{(N)}-B{R^{-1}}B^T((K+\bar{K})x^{(N)}+K_0x_0)\cr
&+(F+\bar{B}_1{P}_0)x_0\big]dt+\frac{1}{N}\sum_{i=1}^ND dW_i\Big\}\\
=&-\Big\{(A_0+B_0P_0)^T(\Lambda_0x_0+\bar{\Lambda}x^{(N)})+(F+B_1P_0)^T((K+\bar{K})x^{(N)}+K_0x_0)\cr
&-\Gamma_1^TQ((I-\Gamma) x^{(N)}-\Gamma_1 x_0)
+P_0^TR_1(P_0x_0+\bar{P}x^{(N)})
\cr
&-P_0^TL^T{R^{-1}}\big[B^T\big((K+\bar{K})x^{(N)}+K_0x_0\big)+ L(P_0{x}_0+\bar{P}{x}^{(N)})\big]\Big\}
+\sum_{j=0}^Nq_0^jdW_j
\end{aligned}$$
which implies $q_0^0=\bar{\Lambda}_0D_0$, $q_0^j=\frac{1}{N}\bar{\Lambda}D, \ j=1,\cdots,N$, and
\begin{align}\label{eq5c}
&\dot{\Lambda}_0 +\Lambda_0(A_0+B_0P_0)+(A_0+B_0P_0)^T\Lambda_0-\bar{\Lambda}BR^{-1}B^TK_0+\bar{\Lambda}(F+\bar{B}_1P_0)\cr
&+(F+\bar{B}_1P_0)^TK_0+P_0^TR_1P_0+\Gamma_1^TQ\Gamma_1-P_0^TL^TR^{-1}LP_0=0,\ \Lambda_0(T)=\bar{\Gamma}_1^TH\bar{\Gamma}_1,\\
\label{eq6c}
&\dot{\bar{\Lambda}}+\bar{\Lambda}(A+G-\bar{B}_1\bar{P})+(A_0+B_0P_0)^T\bar{\Lambda}-\bar{\Lambda}BR^{-1}B^T(K+\bar{K})+\Lambda_0(G_0+B_0\bar{P})\cr&+(F+\bar{B}_1P_0)^T(K+\bar{K})
-Q_{\Gamma_1}^T+P_0^TR_1\bar{P}-P_0^TL^TR^{-1}L\bar{P}=0,\  \bar{\Lambda}(T)=-H_{\bar{\Gamma}_1}^T.
\end{align}
From (\ref{eq5b})-(\ref{eq6c}), it can be verified that $\bar{\Lambda}=K_0^T$, and $K,K,\Lambda_0$ are symmetric matrices.

\begin{theorem}\label{thm2.2}
  Assume that (A1), (A2) and (A4) hold. Problem (P3) admits a solution
  $$\breve{u}_i=-R^{-1}B^T(K{x}_i+\bar{K}{x}^{(N)}+K_0x_0)-R^{-1}L(P_0x_0+\bar{P}x^{(N)}),\ i=1,\cdots,N.$$
\end{theorem}
Let $N\to\infty$. From (\ref{eq63}),
We may approximate $x^{(N)}$ by $\bar{x}$, which satisfies
\begin{equation}\label{eq43}
  d\bar{x}= [(A+G+\bar{B}_1\bar{P}-BR^{-1}B^T(K+\bar{K}))\bar{x}+(F+\bar{B}_1P_0-BR^{-1}B^TK_0)x_0]dt.
\end{equation}
Based on Theorem \ref{thm2.2}, the decentralized feedback strategies for followers can be constructed as:
\begin{equation}
\label{eq43b}\hat{u}_i=-R^{-1}B^T(K{x}_i+\bar{K}\bar{x}+K_0x_0)-R^{-1}L(P_0x_0+\bar{P}\bar{x}).
\end{equation}

\subsection{Optimization for the Leader}

After applying the control laws of followers (\ref{eq43b}),
we have the optimal control problem for the leader.

\textbf{(P4)}: minimize $ {J}_{0}(u_0,\hat{u})$ over $u_0\in \mathcal{U}_d^0$, where
\begin{align}
{J}_0(u_0,\hat{u})=&\mathbb{E}\int_{0}^{T}\big[\|x_0- \Gamma_0\hat{x}^{(N)}\|_{Q_0}^{2}+\|u_0\|_{R_0}^{2}\big]dt
+\mathbb{E}\big[\|x_0(T)- \bar{\Gamma}_0\hat{x}^{(N)}(T)\|_{H_0}^{2}\big],\cr
dx_0=&[A_0x_0+B_0u_0+G_0\hat{x}^{(N)}]dt
+D_0dW_0,\ x_0(0)=\xi_0, 
\\ \label{eq50c}
d\hat{x}_i=&\big\{[A-BR^{-1}B^TK]\hat{x}_i+G\hat{x}^{(N)}-BR^{-1}[B^T\bar{K}+L\bar{P}]\bar{x}\cr
&+[F-BR^{-1}(B^TK_0+LP_0)x_0+B_1u_0]\big\}dt+DdW_i,  \hat{x}_i(0)=\xi_i.
\end{align}

  Since $\{W_i\}$ and $\{x_i(0)\}$ are independent sequences, for the large $N$, it is plausible to replace $\hat{x}^{(N)}$ by $\bar{x}$, which evolves by (\ref{eq43}). 
  In view of (\ref{eq3b}), suppose that the leader has the feedback solution with the  form
$u_0=P_0x_0+\bar{P}\bar{x}.$
This leads to
 a limiting optimal control problem of the leader.



\textbf{(P4$^{\prime}$)}: minimize $ \bar{J}_{0}$ over $P_0,\bar{P}\in C(0,T;\mathbb{R}^{m\times n})$, where
\begin{equation} \label{eq18a}
 \left\{ \begin{aligned}
\bar{J}_0(P_0,\bar{P})=&\mathbb{E}\int_{0}^{T}\big[\|\bar{x}_0- \Gamma_0\bar{x}\|_{Q_0}^{2}+\|P_0x_0+\bar{P}\bar{x}\|_{R_0}^{2}\big]dt+\mathbb{E}\big[\|\bar{x}_0(T)- \bar{\Gamma}_0\bar{x}(T)\|_{H_0}^{2}\big],
\\
d\bar{x}_0=&\big[A_0\bar{x}_0+B_0(P_0\bar{x}_0+\bar{P}\bar{x})+G_0\bar{x}\big]dt
+D_0dW_0,\ x_0(0)=\xi_0, 
\cr
   d\bar{x}=& \big[(A+G-BR^{-1}B^T(K+\bar{K}))\bar{x}\cr
   &+(F-BR^{-1}B^TK_0)\bar{x}_0+\bar{B_1}(P_0\bar{x}_0+\bar{P}\bar{x}) \big]dt, \ \bar{x}(0)=\bar{\xi}.
\end{aligned}\right.
\end{equation}

Let $\bar{X}_0=\mathbb{E}[\bar{x}_0\bar{x}_0^T]$, $\bar{X}=\mathbb{E}[\bar{x}\bar{x}^T]$ and $Y=\mathbb{E}[\bar{x}\bar{x}_0^T]$.
Then by It\^{o}'s formula, we obtain
\begin{align*}
  \frac{d\bar{X}_0}{dt}=&(A_0+B_0P_0)\bar{X}_0+\bar{X}_0(A_0+B_0P_0)^T+(G_0+B_0\bar{P})Y+Y^T(G_0+B_0\bar{P})^T+D_0D_0^T,\\
   \frac{d\bar{X}}{dt}=&(A+G+\bar{B}_1\bar{P}-BR^{-1}B^T(K+\bar{K}))\bar{X}+\bar{X}(A+G+\bar{B}_1\bar{P}-BR^{-1}B^T(K+\bar{K}))^T\cr
   &+(F+\bar{B}_1{P}_0-BR^{-1}B^TK_0)Y^T+Y(F+\bar{B}_1{P}_0-BR^{-1}B^TK_0)^T,\\
  \frac{dY}{dt} =&(A+G+\bar{B}_1\bar{P}-BR^{-1}B^T(K+\bar{K}))Y+(F+\bar{B}_1{P}_0-BR^{-1}B^TK_0)\bar{X}_0\cr
  &+Y(A_0+B_0P_0)^T+\bar{X}(G_0+B_0\bar{P})^T.
    \end{align*}
    The cost of the leader can be expressed equivalently as
    \begin{align*}
      \bar{J}_0(P_0,\bar{P})=&\int_0^Ttr\big(Q_0\bar{X}_0-Q_0\Gamma_0Y-\Gamma_0^TQ_0Y^T+\Gamma_0^TQ_0\Gamma_0\bar{X}+P_0^TR_0P_0\bar{X}_0\cr
      &+\bar{P}^TR_0P_0Y^T+P_0^TR_0\bar{P}Y+\bar{P}^TR_0\bar{P}\bar{X}\big)dt \cr
      &+tr \big[H_0\bar{X}_0(T)-H_0\bar{\Gamma}_0Y(T)-\bar{\Gamma}_0^TH_0Y^T(T)+\bar{\Gamma}_0^TH_0\bar{\Gamma}_0\bar{X}(T)\big].
    \end{align*}
Denote
$$\bar{A}\stackrel{\Delta}{=}A+G+\bar{B}_1\bar{P}-BR^{-1}B^T(K+\bar{K}), \quad \bar{F}\stackrel{\Delta}{=}F+\bar{B}_1{P}_0-BR^{-1}B^TK_0.$$
    Define the Hamiltonian function of the leader as follow:
     \begin{align*}
     &\mathcal{H}(P_0,\bar{P},\Psi_1,\Psi_2,\Psi_3)\cr
     =&tr\Big\{Q_0\bar{X}_0-Q_0\Gamma_0Y-\Gamma_0^TQ_0Y^T+\Gamma_0^TQ_0\Gamma_0\bar{X}+P_0^TR_0P_0\bar{X}_0+\bar{P}^TR_0P_0Y^T\cr
      &+P_0^TR_0\bar{P}Y+\bar{P}^TR_0\bar{P}\bar{X}+[(A_0+B_0P_0)\bar{X}_0+\bar{X}_0(A_0+B_0P_0)^T+(G_0+B_0\bar{P})Y\cr&
      +Y^T(G_0+B_0\bar{P})^T+D_0D_0^T]\Psi_1^T+[\bar{A}\bar{X}+\bar{X}\bar{A}^T+\bar{F}Y^T+Y\bar{F}^T]\Psi_2^T\cr
   &+\big[\bar{A}Y+\bar{F}\bar{X}_0+Y(A_0+B_0P_0)^T+\bar{X}(G_0+B_0\bar{P})^T\big]\Psi_3^T\cr
   &+\big[\bar{A}Y+\bar{F}\bar{X}_0+Y(A_0+B_0P_0)^T+\bar{X}(G_0+B_0\bar{P})^T\big]^T\Psi_3\Big\}.
    \end{align*}
By the matrix maximum principle \cite{A68}, we obtain  the following adjoint equations:
  \begin{align}
     \dot{\Psi}_1=&-\frac{\partial \mathcal{H}}{\partial \bar{X}_0}=-[Q_0+P_0^TR_0P_0+(A_0+B_0P_0)^T\Psi_1+\Psi_1^T (A_0+B_0P_0)+\bar{F }^T\Psi_3+\Psi_3^T\bar{F }],\label{eq27}\\
 \dot{\Psi}_2=&-\frac{\partial \mathcal{H}}{\partial \bar X}=-[\Gamma_0^TQ_0\Gamma_0+\bar{P}^TR_0\bar{P}+\bar{A}^T\Psi_2+\Psi_2\bar{A}\cr
 &+\Psi_3(G_0+B_0\bar{P})+(G_0+B_0\bar{P})^T\Psi_3^T],\label{eq28}\\
  \dot{\Psi}_3=&-\frac{1}{2}\frac{\partial \mathcal{H}}{\partial Y}  =-[\bar{P}^TR_0P_0-\Gamma_0^TQ_0+(G_0+B_0\bar{P})\Psi_1+\Psi_2\bar{F}+\bar{A}^T\Psi_3+\Psi_3(A_0+B_0P_0)]\label{eq29}
  \end{align}
  with the stationarity conditions
    \begin{align}
   0=&\frac{\partial \mathcal{H}}{\partial P_0}=
   2\big(R_0P_0\bar{X}_0+R_0\bar{P}Y+B^T_0\Psi_1\bar{X}_0+ B_0^T\Psi_3^TY\big),\label{eq25}\\
   0=&\frac{\partial \mathcal{H}}{\partial \bar{P}}=
   2\big(R_0P_0Y^T+ R_0\bar{P}\bar{X}+B_0^T\Psi_1Y^T+B_0^T\Psi_3^T\bar{X}\big).\label{eq26}
    \end{align}

  Note that $\Psi_1$ and $\Psi_2$ are symmetric. From (\ref{eq25}) and (\ref{eq26}), we have
\begin{equation}\label{eq82}
  \left\{\begin{aligned}
   P_0&=-R_0^{-1}B^T_0\Psi_1,\\
    \bar{P}&=-R_0^{-1}B^T_0\Psi_3^T.
  \end{aligned}\right.
  \end{equation}
  By applying this into (\ref{eq27})-(\ref{eq29}), we have
   \begin{align*}
      \dot{\Psi}_1=&-[A^T_0\Psi_1+\Psi_1 A_0-\Psi_1B_0R_0^{-1}B^T_0\Psi_1+\bar{F }^T\Psi_3+\Psi_3^T\bar{F }+Q_0],\ \Psi_1(T)=H_0,  
      \\
 \dot{\Psi}_2=&-[\bar{A}^T\Psi_2+\Psi_2\bar{A}-\Psi_3B_0R_0^{-1}B_0^T\Psi_3^T+\Gamma_0^TQ_0\Gamma_0+G_0^T\Psi_3^T+\Psi_3 G_0],\ \Psi_2(T)=\bar{\Gamma}_0^TH_0\bar{\Gamma}_0, 
 \\
  \dot{\Psi}_3=&-[\bar{A}^T\Psi_3+\Psi_3A_0-\Psi_3B_0R_0^{-1}B_0^T\Psi_1+G_0\Psi_1+\Psi_2\bar{F}-\Gamma_0^TQ_0], \ \Psi_3(T)=-\bar{\Gamma}_0^TH_0.
  \end{align*}
  Thus,
the following equation system is obtained as
\begin{equation}\label{eq32a}
\left\{\begin{aligned}
 &\dot{K}+ A^TK+KA-KBR^{-1}B^TK+Q=0,\ K(T)=H,\\
&\dot{\bar{K}}+(A+G+\bar{B}_1\bar{P})^T\bar{K}+\bar{K}(A+G+\bar{B}_1\bar{P})-KBR^{-1}B^T\bar{K}-\bar{K}BR^{-1}B^TK\cr
&-\bar{K}BR^{-1}B^T\bar{K}+(G+\bar{B}_1\bar{P})^TK+K(G+\bar{B}_1\bar{P})+K_0(G_0+B_0\bar{P})
\cr
&+(G_0+B_0\bar{P})^T\bar{\Lambda}+\bar{P}^TR_1\bar{P}-Q_{\Gamma}-\bar{P}^TL^TR^{-1}L\bar{P}=0,\ \bar{K}(T)=-H_{\bar{\Gamma}},\\
&\dot{K}_0+(A+G+\bar{B}_1\bar{P})^TK_0+K_0(A_0+B_0P_0)-(K+\bar{K})BR^{-1}B^TK_0+(G_0+B_0\bar{P})^T\Lambda_0\cr&
+(K+\bar{K})(F+\bar{B}_1P_0)+(\Gamma-I)^TQ\Gamma_1+\bar{P}^TR_1P_0-\bar{P}^TL^TR^{-1}LP_0=0,\  {K}_0(T)=-H_{\bar{\Gamma}_1}\\
&\dot{\Lambda}_0 +\Lambda_0(A_0+B_0P_0)+(A_0+B_0P_0)^T\Lambda_0-\bar{\Lambda}BR^{-1}B^TK_0+\bar{\Lambda}(F+\bar{B}_1P_0)\cr
&+(F+\bar{B}_1P_0)^TK_0+P_0^TR_1P_0+\Gamma_1^TQ\Gamma_1-P_0^TL^TR^{-1}LP_0=0,\ \Lambda_0(T)=\bar{\Gamma}_1^TH_0\bar{\Gamma}_1,\\
&\dot{\bar{\Lambda}}+\bar{\Lambda}(A+G-\bar{B}_1\bar{P})+(A_0+B_0P_0)^T\bar{\Lambda}-\bar{\Lambda}BR^{-1}B^T(K+\bar{K})
+\Lambda_0(G_0+B_0\bar{P})\cr&+(F+\bar{B}_1P_0)^T(K+\bar{K})
+\Gamma_1^TQ(\Gamma-I)+P_0^TR_1\bar{P}-P_0^TL^TR^{-1}L\bar{P}=0,\  \bar{\Lambda}(T)=-H_{\bar{\Gamma}_1}\\
 &\dot{\Psi}_1+A^T\Psi_1+\Psi_1 A-\Psi_1B_0R_0^{-1}B^T_0\Psi_1+\bar{F }^T\Psi_3+\Psi_3^T\bar{F }+Q_0=0,\ \Psi_1(T)=H_0, \\
& \dot{\Psi}_2+\bar{A}^T\Psi_2+\Psi_2\bar{A}-\Psi_3B_0R_0^{-1}B_0^T\Psi_3^T+\Gamma_0^TQ_0\Gamma_0+G_0^T\Psi_3^T+\Psi_3 G_0=0,\ \Psi_2(T)=\bar{\Gamma}_0^TH_0\bar{\Gamma}_0,
 \\
&  \dot{\Psi}_3+\bar{A}^T\Psi_3+\Psi_3A_0-\Psi_3B_0R_0^{-1}B_0^T\Psi_1+G_0\Psi_1+\Psi_2\bar{F}-\Gamma_0^TQ_0=0,\ \Psi_3(T)=-\bar{\Gamma}_0^TH_0,
\end{aligned}\right.
\end{equation}
where $\bar{A}=A+G+\bar{B}_1\bar{P}-BR^{-1}B^T(K+\bar{K})$, and $\bar{F}=F+\bar{B}_1{P}_0-BR^{-1}B^TK_0$.
Based on the above discussions, we may construct the following feedback strategies:
\begin{equation}\label{eq33a}
\left\{\begin{aligned}
\hat{u}_0=&-R_0^{-1}B^T_0(\Psi_1x_0+\Psi_3\bar{x}),\cr
\hat{u}_i=&-R^{-1}B^T(K{x}_i+K_0x_0+\bar{K}\bar{x}) -R^{-1}L(P_0x_0+\bar{P}\bar{x})), \ i=1,\cdots,N,
\end{aligned}\right.
\end{equation}
where $\Psi_1,\Psi_3, 
K,K_0$, $\bar{K}$ are determined by (\ref{eq32a}), $P_0$ and $\bar{P}$ are given by (\ref{eq82}),
and $\bar{x}$ satisfies (\ref{eq43}).

%
%
%

\begin{theorem}\label{thm4.3}

Assume that (A1), (A2) and (A4)  hold, and (\ref{eq32a}) admits a solution. Then the strategy (\ref{eq33a}) is a feedback  ($\epsilon_1,\epsilon_2$)-Stackelberg equilibrium, where $\epsilon_1=\epsilon_2=O(\frac{1}{\sqrt{N}})$. The corresponding social cost of followers is given by
\begin{align}\label{eq77}
  J_{\rm soc}(\hat{u},\hat{u}_0)=&\sum_{i=1}^N\mathbb{E}[\|\xi_i\|^2_{K(0)}]+N\mathbb{E}[\|\xi^{(N)}\|^2_{\bar{K}(0)} +\|\xi_0\|^2_{\Lambda_0}+2\xi_0^T{\bar{\Lambda}}\xi^{(N)}]\cr
  &+N\big(D^T{K}D+D_0^T\Lambda_0D_0\big)+D^T\bar{K}D+N\epsilon_1,
  \end{align}
and the asymptotic cost of the leader is
\begin{align}\label{eq78}
\lim_{N\to\infty}J_{0}(\hat{u},\hat{u}_0)=
\mathbb{E}[\xi_0^T\Psi_1(0)\xi_0+\bar{\xi}^T\Psi_2(0)\bar{\xi}+2\bar{\xi}^T\Psi_3(0)\xi_0]+\mathbb{E}\int_0^T\big(D_0^T\Psi_1D_0\big)dt,
 \end{align}
where
$$\begin{aligned}
  \epsilon_1=&\mathbb{E}\int_0^T\Big\{\|(B^T\bar{K}+L\bar{P})(\hat{x}^{(N)}-\bar{x})\|^2_{R^{-1}}-2(\hat{x}^{(N)}-\bar{x})^T
  \bar{P}^T\big[L^T\hat{u}^{(N)}\cr
&+R_1\big(P_0\hat{x}_0+\frac{1}{2}\bar{P}\bar{x}+\frac{1}{2}\bar{P}\hat{x}^{(N)}\big)+B_0^T(\Lambda_0\hat{x}_0+\bar{\Lambda}\hat{x}^{(N)})
+B_1^T(K+\bar{K})\hat{x}^{(N)}+B^TK_0\hat{x}_0\big]\Big\}dt.
\end{aligned}$$
\end{theorem}

\emph{Proof.} See Appendix \ref{app4}. \hfill{$\Box$}

\section{Simulation}

In this section, we give a numerical example to compare the performances of the open-loop and feedback solutions.
The simulation parameters are listed in Table \ref{tab21}. The step size of iteration is selected as 0.001. 
 The initial distributions of states for the leader and followers satisfy  normal distributions $N(10,2)$ and $N(5,1)$, respectively.

\begin{table}[ht]
	\centering
	\caption{Simulation parameters}
	\label{tab21}       
	\begin{tabular}{cccccccccccccccc}
		\hline\noalign{\smallskip}
	& &	$A_0$&$B_0$&$G_0$&$D_0$&$\Gamma_0$&$Q_0$&$R_0$\\
		\noalign{\smallskip}\hline\noalign{\smallskip}
	& &	$-1$&$1$&$0.1$&$1$&$1$&$1$&$1$ \\
		\noalign{\smallskip}\hline\noalign{\smallskip}
		$A$&$B$&$G$&$F$&$B_1$&$D$&$\Gamma$&$\Gamma_1$&$Q$&$R$&$L$&$R_1$\\
		\noalign{\smallskip}\hline\noalign{\smallskip}
		$-1$&$1$&$0.1$&$1$&$1$&$1$&$1$&$1$&$1$&$2$&$2$&$1$\\
		\noalign{\smallskip}\hline\noalign{\smallskip}
	\end{tabular}
\end{table}

\begin{figure}[h]
	\includegraphics[scale=0.6]{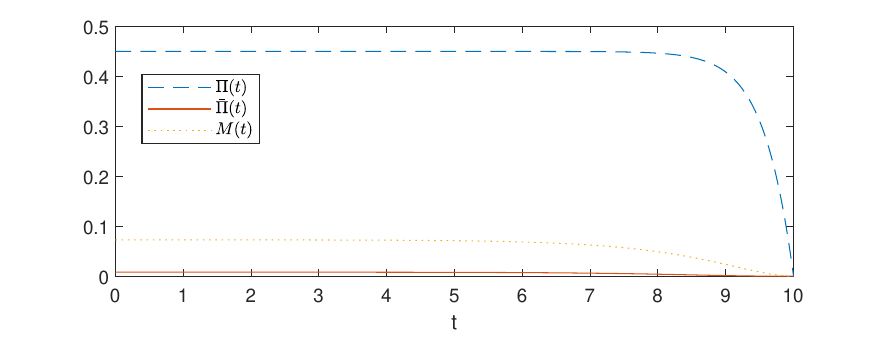}
	\centering
	\caption{ The curves of $\Pi$, $\bar{\Pi}$ and $M$. }
	\label{f11}
\end{figure}

\begin{figure}[h]
	\includegraphics[scale=0.6]{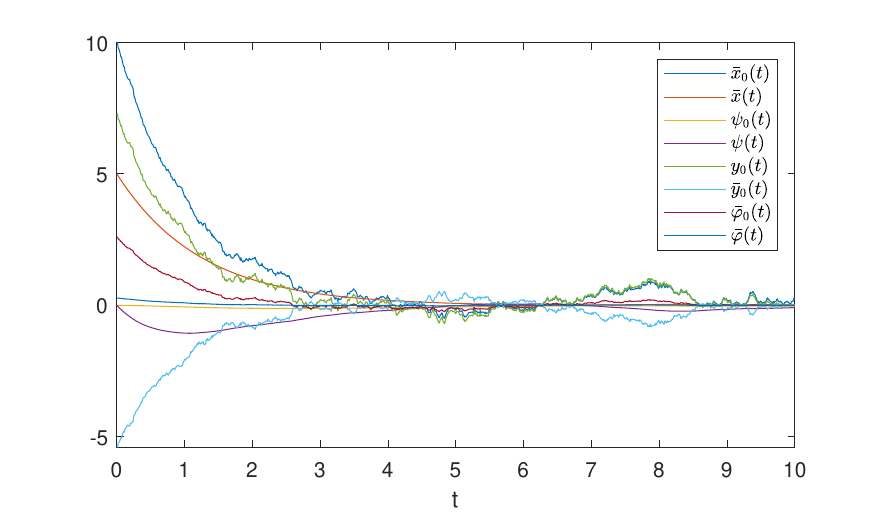}
	\centering
	\caption{ The curves of $X$ and $Y$. }
	\label{f12}
\end{figure}

\begin{figure}[h]
	\includegraphics[scale=0.6]{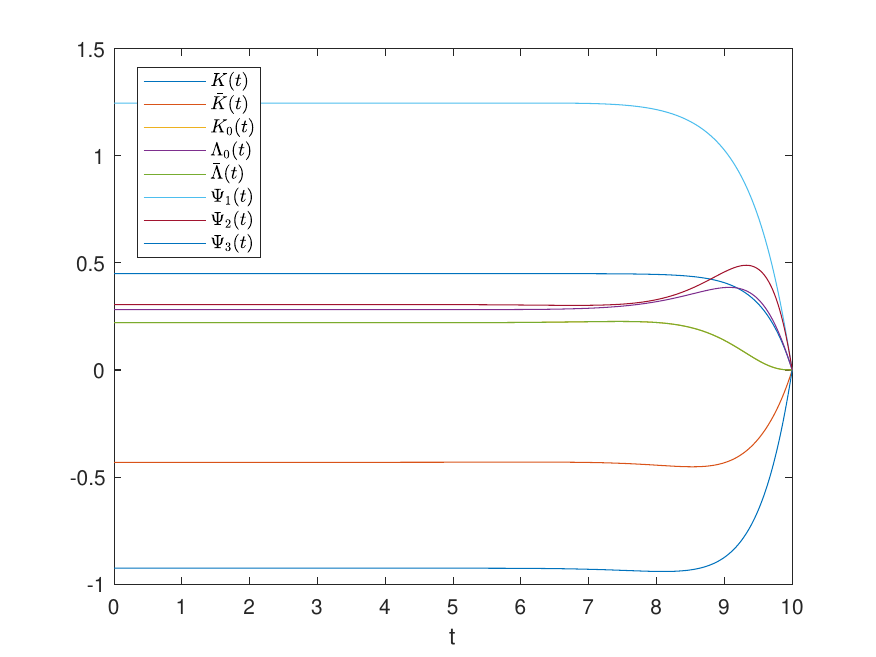}
	\centering
	\caption{ The curves of $K,\bar{K},K_0,$ $\Lambda_0$, $\bar{\Lambda}$, $\Psi_1,\Psi_2$ and $\Psi_3$. }
	\label{f13}
\end{figure}

\begin{figure}[h]
	\includegraphics[scale=0.8]{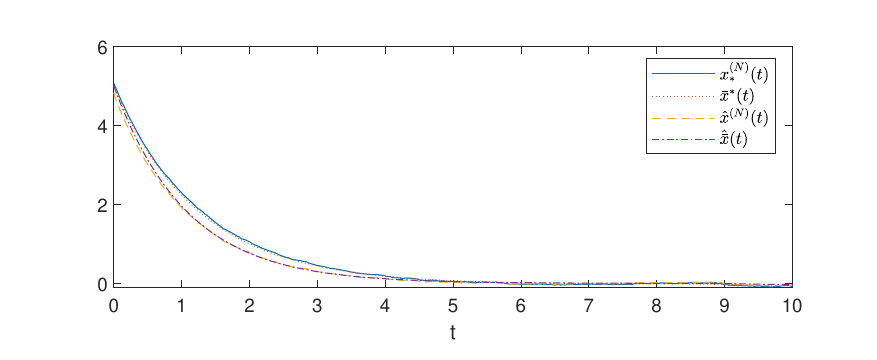}
	\centering
	\caption{  State averages and MF effects of followers under  open-loop and feedback controls. }
	\label{f1}
\end{figure}

\begin{figure}[h]
	\includegraphics[scale=0.8]{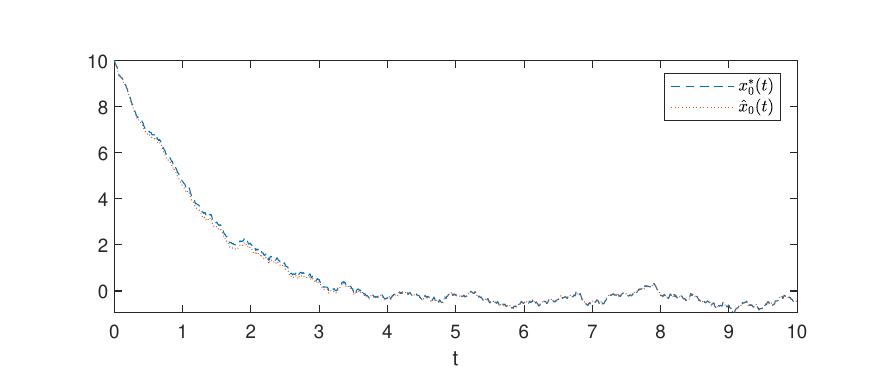}
	\centering
	\caption{State  trajectories of the leader under  open-loop and feedback controls.}
	\label{f2}
\end{figure}

We first  examine the effectiveness of the open-loop and feedback solutions.  Consider a large population system with $1$ leader and $20$ followers.
The decentralized open-loop control \eqref{eq27c} is given  by solving  \eqref{eq6},  \eqref{eq6-b}, \eqref{eq5} and (\ref{eq18}).
Specifically, 
the curves of $\Pi$, $\bar{\Pi}$ and $M$ are shown in Fig. \ref{f11} by virtue of \eqref{eq6}, \eqref{eq6-b} and \eqref{eq5}.  From the Riccati equation \eqref{eq18}, the curves of $X$ and  $Y$ are shown in Fig. \ref{f12}.
The decentralized feedback strategy \eqref{eq33a} is obtained by solving \eqref{eq32a}, and the curves of $K,\bar{K},K_0,$ $\Lambda_0$, $\bar{\Lambda}$, $\Psi_1,\Psi_2$ and $\Psi_3$ are shown in Fig. \ref{f13}. Fig. \ref{f1} gives the trajectories of followers' state averages and MF effects under  open-loop and feedback solutions. Fig. \ref{f2} shows the state  trajectories of the leader under the two solutions. It can be seen  that  state averages approximate MF effects well, and the state average under open-loop control
is larger than the one under feedback control. 

\begin{figure}[h]
	\includegraphics[scale=0.8]{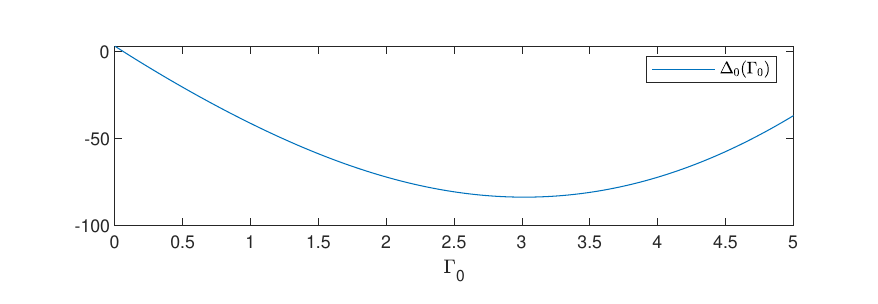}
	\centering
	\caption{ The performance difference $\Delta_0$ for the leader under open-loop and feedback solutions. }
	\label{f3}
\end{figure}

\begin{figure}[h]
	\includegraphics[scale=0.8]{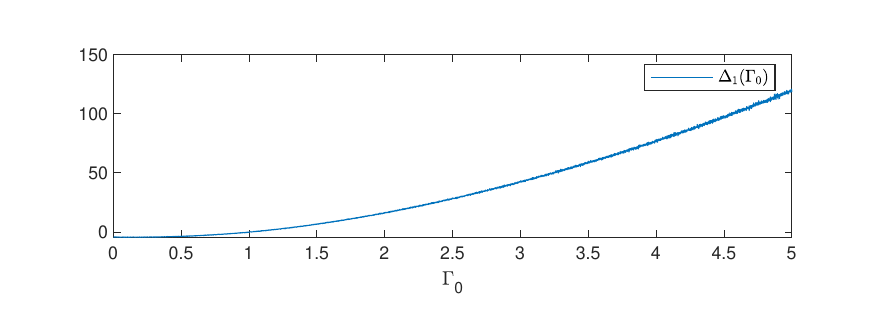}
	\centering
	\caption{The  performance difference $\Delta_1$ for followers under open-loop and feedback solutions.}
	\label{f4}
\end{figure}

Next, we compare the performance of the open-loop and feedback solutions. 
We mainly focus on the influence of the leader' parameter $\Gamma_0$. Denote
\begin{align*}
	\Delta_0(\Gamma_0)\stackrel{\Delta}{=}&J_0(u_0^*,u^*)-J_0(\hat{u}_0,\hat{u}),\\
	\Delta_1(\Gamma_0)\stackrel{\Delta}{=}&\frac{1}{N}[J_{\rm soc}^{ (N)}(u_0^*,u^*)-J_{\rm soc}^{ (N)}(\hat{u}_0,\hat{u})].
\end{align*}
Start with $\Gamma_0=0$, and increase $\Gamma_0$ by $0.001$ at each step until $\Gamma_0=5$. We calculate the differences $\Delta_0$ and $\Delta_1$ to compare the  performance of the two solutions with respect to different $\Gamma_0$. For each $\Gamma_0$, we compute $200$ times and take the average. The trajectories of $\Delta_0(\Gamma_0)$ and $\Delta_1(\Gamma_0)$ are plotted in Figs. \ref{f3} and \ref{f4}, respectively. The figures show that in  above parameter setting,  the open-loop
control engenders a lower cost 
 for the leader than the feedback control. However, the opposite is
true for followers, and the feedback solution generally outperforms the open-loop
solution.


\section{Concluding remarks}
This paper studies open-loop and feedback solutions of LQ-MF  games with a leader and a large number
of followers. By
variational analysis with MF  approximations, we obtain a set of  open-loop
controls of players in terms of solutions to MF  FBSDEs. By applying the matrix maximum principle, a
set of decentralized feedback strategies is constructed for all the players. For both solutions, the corresponding optimal costs of all players are explicitly given by
virtue of the solutions to two Riccati equations, respectively.
In further works, it is very interesting to investigate model-free MF  Stackelberg games by reinforcement learning.

\appendices

\section{Proofs of Theorems \ref{thm1} and \ref{thm3.3}}\label{app1}
\def\theequation{A.\arabic{equation}}
\setcounter{equation}{0}
\emph{Proof of Theorem \ref{thm1}}.
Suppose that $\check{u}_i=-R^{-1}(B^T\check{p}_i+ L\check{u}_0),$ where {${(\check{p}_i, \check{q}_{i}^{j}},i,j=0,1,\cdots,N)$ is a set of solutions to the backward equations in (\ref{eq3a}).}
 Denote by $\check{x}_0$, $\check{x}_i$ the state of agent $i$ under the control $\check{u}_i$, $i=1,\cdots,N$. For any $u_i\in L^2_{{\mathcal  F}}(0, T; \mathbb{R}^r) $
 and $\theta\in \mathbb{R}\ (\theta\not=0)$, let $u_i^{\theta}=\check{u}_i+\theta u_i$.
 Denote by $ x_i^{\theta}$, $i=0,1,\cdots,N$ the solution of the following perturbed state equation
$$\left\{ \begin{aligned}
 dx_0^{\theta}=&[A_0x_0^{\theta}+B_0u_0
+\frac{G_0}{N}\sum_{i=1}^Nx^{\theta}_i]dt
+D_0dW_0,
\\
dx_i^{\theta}=&\Big[Ax_i^{\theta}+B(\check{u}_i+\theta u_i)+\frac{G}{N}\sum_{i=1}^Nx^{\theta}_i+Fx^{\theta}_0\Big]dt+DdW_i,\cr
x_0(0)&=\xi_{0},\ x_i(0)=\xi_{i},\ i=1,2,\cdots,N.
 \end{aligned}\right.$$
 Let $z_i=(x_i^{\theta}-\check{x}_i)/\theta$, $i=0,1,2,\cdots,N$. 
It can be verified that
 $z_0$ and $z_i$ satisfy
 $$\left\{ \begin{aligned}
 dz_0=&\big[A_0z_0+{G_0}z^{(N)}\big]dt,
\\
dz_i=&\big[Az_i+B u_i+{G}z^{(N)}+Fz_0\big]dt,\cr
z_0(0)&= 0,\ z_i(0)=0,\ i=1,2,\cdots,N.
 \end{aligned}\right.$$
Then by It\^{o}'s formula, we have
{\begin{align}\label{eq-6a}
&\mathbb{E}[\langle- H_{\bar{\Gamma}_1} ^Tx^{(N)}(T)+\bar{\Gamma}_1^TH\bar{\Gamma}_1x_0(T),z_0(T)\rangle ]
=\mathbb{E}[\langle \check{p}_0(T),z_0(T)\rangle-\langle \check{p}_0(0),z_0(0)\rangle]\cr
=&\mathbb{E}\!\int_0^T\! \Big\{\big\langle -\big[A^T_0\check{p}_0+F^T\check{p}^{(N)}-\Gamma_1Q\big((I-\Gamma) \check{x}^{(N)}-\Gamma_1 \check{x}_0\big)\big],
z_0\rangle\cr
&+\langle \check{p}_0,A_0z_0+{G_0}z^{(N)}\rangle\Big\} dt\cr
=&\mathbb{E}\!\int_0^T\!  \Big\{\big\langle -\big[F^T\check{p}^{(N)}-\Gamma_1Q\big((I-\Gamma) x^{(N)}-\Gamma_1 x_0\big)\big],
z_0\rangle+\langle \check{p}_0,{G_0}z^{(N)}\rangle \Big\} dt,
\end{align}}
and 
{\begin{align}\label{eq-7b}
&\sum_{i=1}^N\mathbb{E}[\langle Hx_i(T)-H_{\bar{\Gamma}}x^{(N)}(T)-H_{\bar{\Gamma}_1} x_0(T),z_i(T)\rangle]\cr
=&\sum_{i=1}^N\mathbb{E}[\langle \check{p}_i(T),z_i(T)\rangle-\langle \check{p}_i(0),z_i(0)\rangle]\cr
=&\sum_{i=1}^N\mathbb{E}\!\int_0^T\!  \Big\{\big\langle -\big[A^T\check{p}_i+G^T\check{p}^{(N)}+G^T_0\check{p}_0+Q\check{x}_i-Q_{\Gamma}\check{x}^{(N)}\cr
&
+(\Gamma-I)^TQ\Gamma_1\check{x}_0\big],z_i\big\rangle+\langle \check{p}_i,Az_i+Gz^{(N)}+Fz_0+Bu_i\rangle \Big\} dt\cr
= &\sum_{i=1}^N \mathbb{E}\int_0^T \Big\{\big\langle -\big[G^T_0\check{p}_0+Q\check{x}_i-Q_{\Gamma}\check{x}^{(N)}
+(\Gamma-I)^TQ\Gamma_1\check{x}_0\big],z_i\big\rangle\cr
&+\langle F^T\check{p}^{(N)},z_0\rangle+\langle B^T\check{p}_i,u_i\rangle\Big\} dt.
\end{align}}
From (\ref{eq3}),
\begin{equation}\label{eq5a}
\begin{aligned}
 &{J}_{\rm soc}(\check{u}+\theta u)-{J}_{\rm soc}(\check{u})=2\theta I_1+{\theta^2}I_2
\end{aligned}
\end{equation}
where $\check{u}=(\check{u}_1,\cdots,\check{u}_N)$, and
\begin{align*}
I_1\stackrel{\Delta}{=}&\sum_{i=1}^N\mathbb{E}\int_0^T \big[\big\langle Q\big(\check{x}_i-(\Gamma\check{x}^{(N)}+\Gamma_1\check{x}_0)\big),
 z_i
-\Gamma z^{(N)}-\Gamma_1z_0\big\rangle
 +
\langle R \check{u}_i+L\check{u}_0,u_i\rangle \big]dt\cr
&+\sum_{i=1}^N\mathbb{E}[\langle H\big(\check{x}_i(T)-(\bar{\Gamma}\check{x}^{(N)}(T)+\bar{\Gamma}_1\check{x}_0(T))\big),
 z_i(T)
-\bar{\Gamma} z^{(N)}(T)-\bar{\Gamma}_1z_0(T)\big\rangle],\cr
I_2\stackrel{\Delta}{=}&\sum_{i=1}^N\mathbb{E}\int_0^T\! \big[\big\|z_i
  \! -\!\Gamma z^{(N)}-\Gamma_1z_0\big\|^2_{Q}
+\|u_i\|^2_{R}\big]dt+\sum_{i=1}^N\mathbb{E}[\big\|z_i(T)
  \! -\bar{\Gamma} z^{(N)}(T)-\bar{\Gamma}_1z_0(T)\big\|^2_{H}].
\end{align*}
Note that
\begin{align*}
  &\sum_{i=1}^N\mathbb{E}\Big[\int_0^T \big\langle Q\big(\check{x}_i-(\Gamma\check{x}^{(N)}+\Gamma_1\check{x}_0)\big),\Gamma z^{(N)}\big\rangle dt\cr
  &+\big\langle H\big(\check{x}_i(T)-(\bar{\Gamma}\check{x}^{(N)}(T)+\bar{\Gamma}_1\check{x}_0(T))\big),\bar{\Gamma } z^{(N)}(T)\big\rangle \Big]\cr
  =&
\mathbb{ E}\!\int_0^T\! \Big\langle \Gamma^TQ  \sum_{i=1}^N\big(\check{x}_i\!-\!(\Gamma\check{x}^{(N)}\!+\Gamma_1\check{x}_0)\big),\frac{1}{N}  \sum_{j=1}^Nz_j\Big\rangle  dt\cr
&+\mathbb{ E}\Big[\big\langle\bar{\Gamma }^T H\sum_{i=1}^N\big(\check{x}_i(T)-(\bar{\Gamma}\check{x}^{(N)}(T)+\bar{\Gamma}_1\check{x}_0(T))\big), \frac{1}{N}  \sum_{j=1}^Nz_j(T)\big\rangle \Big]\cr
 =& \sum_{j=1}^N \mathbb{E}\Big[\int_0^T \big\langle  {\Gamma^TQ} \big((I-\Gamma)\check{x}^{(N)}-\Gamma_1\check{x}_0\big), z_j\big\rangle  dt+\big\langle \bar{\Gamma }^TH\big((I-\bar{\Gamma})\check{x}^{(N)}(T)-\bar{\Gamma}_1\check{x}_0(T)\big), z_j(T)\big\rangle\Big],
\end{align*}
and
$$\sum_{i=1}^N \mathbb{E}\int_0^T \langle \check{p}_0,{G_0}z^{(N)}\rangle dt=\sum_{i=1}^N \mathbb{E}\int_0^T \langle G_0^T\check{p}_0,z_i\rangle dt.$$
In views of (\ref{eq-6a}) and (\ref{eq-7b}), by direct computations we obtain that
{\begin{align}\label{eq10b}
I_1=&\sum_{i=1}^N\mathbb{E}\int_0^T\big[\big\langle Q\big(\check{x}_i-(\Gamma\check{x}^{(N)}+\Gamma_1\check{x}_0)\big)-\Gamma^TQ \big((I-\Gamma)\check{x}^{(N)}-\Gamma_1\check{x}_0\big),\ z_i\big\rangle\cr
&-\big\langle \Gamma_1^TQ\big(\check{x}_i-(\Gamma\check{x}^{(N)}+\Gamma_1\check{x}_0)\big),z_0\big\rangle +
\langle R\check{u}_i+ L\check{u}_0+B^T\check{p}_i,u_i\rangle \big]dt\cr
&+\sum_{i=1}^N \mathbb{E}\int_0^T  \big\{\langle -\big[F^T\check{p}^{(N)}-\Gamma_1^TQ((I-\Gamma) \check{x}^{(N)}-\Gamma_1 \check{x}_0)\big],
z_0\rangle \cr
&-\langle Q\check{x}_i-Q_{\Gamma}\check{x}^{(N)}
+(\Gamma-I)^TQ\Gamma_1\check{x}_0,z_i\big\rangle +\langle F^T\check{p}^{(N)},z_0\rangle \big\}dt\cr
=&\sum_{i=1}^N\mathbb{E}\int_0^T\big[\big\langle Q\check{x}_i-Q_{\Gamma}\check{x}^{(N)}+(\Gamma-I)^TQ\Gamma_1\check{x}_0), z_i\big\rangle\cr
&-\big\langle \Gamma_1^TQ\big((I-\Gamma)\check{x}^{(N)}-\Gamma_1\check{x}_0\big),z_0\big\rangle +
\langle R\check{u}_i+ L\check{u}_0+ B^T\check{p}_i,u_i\rangle \big]dt\cr
&+\sum_{i=1}^N \mathbb{E}\int_0^T  \big\{\big\langle \big[\Gamma_1^TQ((I-\Gamma) \check{x}^{(N)}-\Gamma_1 \check{x}_0)\big],
z_0\big\rangle \cr
&-\langle Q\check{x}_i-Q_{\Gamma}\check{x}^{(N)}
+(\Gamma-I)^TQ\Gamma_1\check{x}_0,z_i\big\rangle \big\}dt\cr
=&\sum_{i=1}^N\mathbb{E}\int_0^T\big\langle R\check{u}_i+ L\check{u}_0+B^T\check{p}_i,u_i\big\rangle dt.
\end{align}}
Since $Q\geq0$ and $R>0$, we have $I_2\geq0$ and Problem (P1) admit an optimal control \cite{YZ99}. From (\ref{eq5a}), $\check{u}$ is a minimizer to Problem (P1) if and only if
 $I_1=0 $, which is equivalent to
$\check{u}_i=-{R^{-1}}(B^T\check{p}_i+ L\check{u}_0).$
Thus, we have the optimality system (\ref{eq3}).
This implies that (\ref{eq3}) admits a solution $(\check{x}_i,\check{p}_i,\check{q}_{i}^{j}, i,j=1,\cdots,N)$. \hfill{$\Box$}

\emph{Proof of Theorem \ref{thm3.3}.} Suppose $\{u_0^*\}$ is the optimal control of Problem (P2$^\prime$).
$\bar{x}_0^*$ and $\bar{x}^*$ are the corresponding states in (\ref{eq8b}) and (\ref{eq8c})  under the control $\{{u}_0^*\}$.
 For $i=1,2,\cdots,N,$ denote $\delta \bar{x}_0=\bar{x}_0-\bar{x}_0^*$ the increment of $\bar{x}_0$ along with the variation $\delta u_0=u_0-u_0^*.$ Similarly, $\delta \bar{x}=\bar{x}-\bar{x}^*$, $\delta \bar{\varphi_0}=\bar{\varphi}_0(u_0)-\bar{\varphi}_0^*(u_0^*)$ and $\delta \bar{\varphi}=\bar{\varphi}(u_0)-\bar{\varphi}^*(u_0^*)$.
 Define the variation of $\bar{J}_{0}$ on ${u}^*$:
$$\delta \bar{J}_{0}(\delta u_0 ,u^* )= \bar{J}_{0}(u_0,{u}^*)- \bar{J}_{0}(u_0^*,{u}^*(u_0^*))+o(||\delta u_{0}||_{L^{2}}).$$ Then,
we have
\begin{align*}
  d\delta \bar{x}_0=&(A_0\delta \bar{x}_0+G_0\delta\bar{x}+B_0 \delta u_0)dt,\  \delta x_0(0)=0,\\
   d\delta \bar{x}=&[(A+G-BR^{-1}B^T\bar{\Pi} )\delta \bar{x}-BR^{-1}B^T\delta\bar{\varphi}\cr
   &+(F-BR^{-1}B^TM) \delta \bar{x}_0+\bar{B}_1\delta u_0]dt,  \
   \delta \bar{x}(0)=0,\\
   d\delta{\bar{\varphi}_0}=&-\big[A_0^T\delta\bar{\varphi}_0+(F-BR^{-1}B^T\Pi_0)^T\delta\bar{\varphi}+(\Pi_0\bar{B}_1+M_0B_0)\delta u_0\big]dt+\delta\bar{q}_0^0dW_0,\  \varphi_0(T)=0,\cr
   d\delta \bar{\varphi}=&-\big[(A+G-BR^{-1}B^T\bar{\Pi})^T\delta \bar{\varphi}+G^T_0 \delta \bar{\varphi}_0
\cr&+(\bar{\Pi}\bar{B}_1+MB_0)\delta u_0\big]dt+\delta\bar{q}_i^0dW_0, \  \delta \bar{\varphi}(T)=0.
\end{align*}
Let $y_0,\bar{y}$,
 $\psi_0$, and $\psi$ satisfy (\ref{eq16}).
By It\^{o}'s formula, we obtain
\begin{align}\label{eq12}
&\mathbb{E}[\langle H_0(\bar{x}_0^*(T)-\bar{\Gamma}_0\bar{x}^*(T)),\delta \bar{x}_0(T)\rangle]=\mathbb{E}[\langle y_0(T),\delta \bar{x}_0(T)\rangle-\langle y_0(0),\delta \bar{x}_0(0)\rangle]\cr
=&\mathbb{E}\int_0^T\big[-\langle (F-BR^{-1}B^TM)^T\bar{y}+Q_0(\bar{x}_0^*-\Gamma_0\bar{x}^*),\delta \bar{x}_0\rangle
+\langle G_0^Ty_0,\delta \bar{x}\rangle+\langle B_0^Ty_0,\delta u_0\rangle\big]dt,
\end{align}
and
\begin{align}\label{eq13}
&\mathbb{E}[\langle -\bar{\Gamma}_0^TH_0(\bar{x}_0^*(T)-\bar{\Gamma}_0\bar{x}^*(T)),\delta \bar{x}(T)\rangle]
 =\mathbb{E}[\langle \bar{y}(T),\delta \bar{x}(T)\rangle-\langle \bar{y}(0),\delta \bar{x}(0)\rangle]\cr
 =&\mathbb{E}\int_0^T[\langle \Gamma_0^TQ_0 (\bar{x}_0^*-\Gamma_0\bar{x}^*)-G_0^Ty_0,\delta \bar{x}\rangle
-\langle BR^{-1}B^T \bar{y},\delta \bar{\varphi}\rangle]dt\cr
&+\mathbb{E}\int_0^T[\langle (F-BR^{-1}B^TM)^T\bar{y},\delta \bar{x}_0\rangle+\langle \bar{B}_1^T\bar{y},\delta u_0\rangle]dt,
\end{align}
Similarly, we have
\begin{align}\label{eq14}
 0= &\mathbb{E}[\langle\delta\bar{\varphi}_0(T),\psi_0(T)\rangle -\langle\delta\bar{\varphi}_0(0),\psi_0(0)\rangle]\cr
=& \mathbb{E}\int_0^T[\langle G_0\psi, \delta\bar{\varphi}_0\rangle-\langle(F-BR^{-1}B^T\Pi_0)\psi_0, \delta\bar{\varphi}\rangle \cr
&-\langle({\Pi}_0\bar{B}_1+M_0B_0)^T\psi_0,\delta u_0\rangle+\langle \delta \bar{q}_0^0,v_0\rangle] dt,
\end{align}
and
\begin{align}\label{eq14b}
 0= &\mathbb{E}[\langle\delta\bar{\varphi}(T),\psi(T)\rangle -\langle\delta\bar{\varphi}(0),\psi(0)\rangle]\cr
=& \mathbb{E}\int_0^T[\langle(F-BR^{-1}B^T\Pi_0)\psi_0-BR^{-1}B^T \bar{y}, \delta\bar{\varphi}\rangle -\langle G_0 \psi,\delta \bar{\varphi}_0\rangle\cr
&-\big\langle(\bar{\Pi}\bar{B}_1+MB_0)^T\psi,\delta u_0\big\rangle+\langle \delta \bar{q}_i^0,v\rangle] dt.
\end{align}
From  (\ref{eq12})-(\ref{eq14b}),
$$\begin{aligned}
\delta{J}_0(u_0,u^*)=&\mathbb{E}\int_0^T[\langle Q_0(\bar{x}_0^*-\Gamma_0\bar{x}^*),\delta \bar{x}_0-\Gamma_0\delta \bar{x}\rangle+\langle u_0^*, R_0\delta u_0\rangle]dt\cr&+\mathbb{E}[\langle H_0(\bar{x}_0^*(T)-\bar{\Gamma}_0\bar{x}^*(T)),\delta \bar{x}_0(T)-\bar{\Gamma}_0\delta \bar{x}(T)\rangle]\cr
=&\mathbb{E}\int_0^T\langle B_0^Ty_0+\bar{B}_1^T\bar{y}+({\Pi}_0\bar{B}_1+M_0B_0)^T\psi_0\cr
  &+(\bar{\Pi}\bar{B}_1+MB_0)^T\psi+Ru_0^*,\delta u_0\rangle dt
=0.
\end{aligned}$$
Note that $Q_0\geq 0$ and $R_0>0$. By the standard variational principle \cite{YZ99}, $u_0^*$ is an optimal control of (P2$^{\prime}$).
\hfill $\Box$

\section{Proof of Theorem \ref{thm3.4}.}\label{app2}
\def\theequation{B.\arabic{equation}}
\setcounter{equation}{0}

To prove Theorem \ref{thm3.4}, we first give a lemma.  Let $x_i^*$ be the realized state under the control $u^*_i, i=0,1,\cdots,N$.
Denote ${x}_*^{(N)}=\frac{1}{N}\sum_i^N{x}_i^*$.
 \begin{lemma}\label{lem3.1}
Under (A1)-(A3), the following hold:
 \begin{equation}\label{eq39}
\sup_{0\leq t\leq T}\mathbb{E}\big[\| {x}_0^*-\bar{x}_0\|^2+\| {x}_*^{(N)}-\bar{x}\|^2\big]dt=O(1/N),
\end{equation}
\begin{equation}\label{eq39b}
\sup_{0\leq t\leq T}\mathbb{E}\big[\| x_i^*-\bar{x}_i\|^2\big]dt=O(1/N).
\end{equation}
 \end{lemma}
 \emph{Proof.} By (\ref{eq8c}) and (\ref{eq40}), it can be verified that
 $\sup_{0\leq t\leq T}\mathbb{E}\big[\|\bar{x}^{(N)}-\bar{x}\|^2\big]dt=O(1/N). $
From (\ref{eq18d}), (\ref{eq23})-(\ref{eq8c}), we have
$$\begin{aligned}
d({x}_0^*-\bar{x}_0)=&[A_0({x}_0^*-\bar{x}_0)+G_0({x}_*^{(N)}-\bar{x})]dt,\ x_0^*(0)-\bar{x}_0(0)=0,\\
  d({x}_*^{(N)}-\bar{x})=&\big[(A+G)
({x}_*^{(N)}-\bar{x})+F({x}_0^*-\bar{x}_0)-BR^{-1}B^T\Pi(\bar{x}^{(N)}-\bar{x})\big]dt\cr
&+\frac{1}{N}\sum_{j=1}^NDdW_j,\ {x}_*^{(N)}(0)-\bar{x}(0)=\xi^{(N)}-\bar{\xi},
\end{aligned} $$
where $\xi^{(N)}=\frac{1}{N}\sum_{j=1}^N\xi_j$.
Let $\mathbb{A}=\left[\begin{array}{cc}
 A_0&G_0\\
 F&A+G-BR^{-1}B^T\Pi
\end{array}\right]$. Let $\Phi$ be the solution of the equation
$\dot{\Phi}=\mathbb{A}\Phi, \ \Phi(0)=I$.
Then we have
\begin{align*}
  \left[\begin{array}
  {c}{x}_0^*(t)-\bar{x}_0(t)\\
{x}_*^{(N)}(t)-\bar{x}(t)
\end{array}\right] \leq &\Phi(t) \left[\begin{array}{c}
  0\\
\xi^{(N)}-\bar{\xi}
\end{array}\right]+ \int_0^t\Phi(t-s)\left[\begin{array}{c}
  0\\
\frac{1}{N}\sum_{j=1}^NDdW_j(s)
\end{array}\right]\cr
&-\int_0^T\Phi(t-s)\left[\begin{array}{c}
  0\\
BR^{-1}B^T\Pi(s)(\bar{x}^{(N)}(s)-\bar{x}(s))
\end{array}\right]ds
\end{align*}
which gives (\ref{eq39}).
Note that $x_i^*-\bar{x}_i$ satisfies
$$d(x_i^*-\bar{x}_i)=[A(x_i^*-\bar{x}_i)+G({x}_*^{(N)}-\bar{x})+F({x}_*^{(N)}-\bar{x})]dt.$$
By (\ref{eq39}), we can obtain (\ref{eq39}). \hfill{$\Box$}

 \emph{Proof of Theorem \ref{thm3.4}.}
 \emph{(For followers).} 
For any $u_i\in \mathcal{U}_c$, let $\tilde{u}_i=u_i-{u}_i^*$, $\tilde{x}_i=x_i-{x}_i^*$, $\tilde{x}_0={x}_0-{x}_0^*$ and $\tilde{x}^{(N)}=\frac{1}{N}\sum_{i=1}^N \tilde{x}_i$. Then by (\ref{eq1}), (\ref{eq18d}) and (\ref{eq19d}),
\begin{equation}\label{eq32-a}
\left\{
\begin{aligned}
d\tilde{x}_0&=(A_0\tilde{x}_0+{G}_0\tilde{x}^{(N)})dt, \  \tilde{x}_0(0)=0,\\
d\tilde{x}_i&=(A\tilde{x}_i+{G}\tilde{x}^{(N)}+B\tilde{u}_i+F\tilde{x}_0)dt, \  \tilde{x}_i(0)=0.
\end{aligned}\right.
\end{equation}
From (\ref{eq3}), we have
$J_{\rm soc}(u_0,u)=\sum_{i=1}^N({J}_i(u_0,{u}^*)+\tilde{J}_i(u_0,\tilde{u})+I_i),$
where
\begin{align*}
\tilde{J}_i(u_0,\tilde{u})\stackrel{\Delta}{=}\mathbb{E}\int_0^T&\big[\|\tilde{x}_i-\Gamma \tilde{x}^{(N)}-\Gamma _1 \tilde{x}_0\|^2_Q+\|\tilde{u}_i\|^2_{R}\big]dt+\mathbb{E}\|\tilde{x}_i(T)-\bar{\Gamma} \tilde{x}^{(N)}(T)-\bar{\Gamma} _1 \tilde{x}_0(T)\|^2_H,\\
\mathcal{I}_i=2\mathbb{E}\int_0^T&\big[\big({x}_i^*-\Gamma {x}_*^{(N)}-\Gamma _1 {x}^*_0\big)^TQ\big(\tilde{x}_i
-\Gamma \tilde{x}^{(N)}-\Gamma _1 \tilde{x}_0\big)+\tilde{u}_i^TLu_0^*+\tilde{u}^T_iR{u}_i^*\big]dt\\
&+\mathbb{E}\big[\big({x}_i^*(T)-\bar{\Gamma} {x}_*^{(N)}(T)-\bar{\Gamma} _1 {x}^*_0(T)\big)^TH\big(\tilde{x}_i(T)
-\bar{\Gamma} \tilde{x}^{(N)}(T)-\bar{\Gamma} _1 \tilde{x}_0(T)\big)\big].
\end{align*}
Let $\bar{p}^{(N)}=\frac{1}{N}\sum_{i=1}^N\bar{p}_i$. Note that $\bar{p}_i-\bar{p}=\Pi(\bar{x}_i-\bar{x})$. 
By (\ref{eq4a})-(\ref{eq11-a}) and It\^{o}'s formula, we obtain
$$  \begin{aligned}
&N\mathbb{E}\big[  \tilde{x}_0^T(T)\big(-H_{\bar{\Gamma}_1} ^T\bar{x}(T)+\bar{\Gamma}_1^TH\bar{\Gamma}_1\bar{x}_0(T)\big)\big]=\sum_{i=1}^N\mathbb{E}\big[  \tilde{x}_0^T(T)\bar{p}_0(T)\big] \\
=&\sum_{i=1}^N\mathbb{E}\int_{0}^{T}\Big\{\tilde{x}_i^T G_0^T\bar{p}_0- \tilde{x}_0^T\big[ F^T\bar{p}^{(N)}+\Gamma^T_1Q((\Gamma-I) \bar{x}+\Gamma_1\bar{x}_0)+F^T\Pi(\bar{x}-\bar{x}^{(N)})\big]\Big\}dt,
\end{aligned}$$
and
$$  \begin{aligned}
&\sum_{i=1}^N\mathbb{E}\big[\tilde{x}_i^T(T)\big(H\bar{x}_i(T)-H_{\bar{\Gamma}}\bar{x}(T)-H_{\bar{\Gamma}_1} \bar{x}_0(T)\big)\big]
=\sum_{i=1}^N\mathbb{E}\big[  \tilde{x}_i^T(T)\bar{p}_i(T))\big] \\
=\ &\sum_{i=1}^N\mathbb{E}\int_{0}^{T}\big[(A\tilde{x}_i+{G}\tilde{x}^{(N)}+B\tilde{u}_i+F\tilde{x}_0)^T\bar{p}_i-\tilde{x}_i^T( A^T\bar{p}_i+G^T\bar{p}+G_0\bar{p}_0\\
&+Q\bar{x}_i-Q_{\Gamma} \bar{x}-(I-\Gamma)^TQ\Gamma_1\bar{x}_0)\big]dt \\
=\ &\sum_{i=1}^N\mathbb{E}\int_{0}^{T}\Big\{\tilde{x}_0^TF^T\bar{p}_i- \tilde{x}_i^T\big[ G_0^T\bar{p}_0+Q\bar{x}_i-Q_{\Gamma} \bar{x}-(I-\Gamma)^TQ\Gamma_1 \bar{x}_0\\
&+G^T\Pi(\bar{x}^{(N)}-\bar{x})\big]-\tilde{u}_i^TB^T\bar{p}_i\Big\}dt,
\end{aligned}$$
which lead to
 \begin{align}\label{eq45}
  &\sum_{i=1}^N\mathbb{E}\big[\tilde{x}_i^T(T)\big(H\bar{x}_i(T)-H_{\bar{\Gamma}}\bar{x}(T)-H_{\bar{\Gamma}_1} \bar{x}_0(T)\big)+\tilde{x}_0^T(T)\big(-H_{\bar{\Gamma}_1} ^T\bar{x}(T)+\bar{\Gamma}_1^TH\bar{\Gamma}_1\bar{x}_0(T)\big)\big]\cr
  =&\sum_{i=1}^N\mathbb{E}\int_{0}^{T}\Big\{- \tilde{x}_0^T\big[ \Gamma_1^TQ\Gamma_1\bar{x}_0-Q_{\Gamma_1} \bar{x})+F^T\Pi(\bar{x}-\bar{x}^{(N)})\big]- \tilde{x}_i^T\big[ Q\bar{x}_i-Q_{\Gamma} \bar{x}\cr 
&-Q_{\Gamma_1} \bar{x}_0+G^T\Pi(\bar{x}^{(N)}-\bar{x})\big]-\tilde{u}_i^TB^T\bar{p}_i\Big\}dt.
\end{align}
Note that $B^T\bar{p}_i=R({u}_i^*+Lu_0^*)$ and
$$\begin{aligned}
  \sum_{i=1}^N \mathcal{I}_i=&  \sum_{i=1}^N2 \mathbb{E}\int_{0}^{T}\Big[\tilde{x}_i^T\big( Q{x}_i^*-Q_{\Gamma} {x}_*^{(N)}-Q_{\Gamma_1} x_0^*\big)+\tilde{x}_0^T(\Gamma_1^TQ\Gamma_1x_0^*-Q_{\Gamma_1}^Tx^{(N)}_*)\cr
 &+\tilde{u}_i^TR({u}_i^*+Lu_0^*)\Big]dt+ \sum_{i=1}^N\mathbb{E}\Big[\tilde{x}_i^T(T)\big(H{x}^*_i(T)-H_{\bar{\Gamma}}{x}_*^{(N)}(T)-H_{\bar{\Gamma}_1} {x}_0^*(T)\big)\cr&
 +\tilde{x}_0^T(T)\big(-H_{\bar{\Gamma}_1}^T{x}_*^{(N)}(T)+\bar{\Gamma}_1^TH\bar{\Gamma}_1x_0^*(T)\big)\Big].
\end{aligned}$$
From Lemma \ref{lem3.1} and (\ref{eq45}), 
one can obtain
\begin{equation*}
\begin{aligned}
\frac{1}{N}\sum_{i=1}^N \mathcal{I}_i
=\ &\frac{1}{N} \sum_{i=1}^N 2\mathbb{E}\int_{0}^{T} \Big\{\tilde{x}_i^T\big[G^T\Pi(\bar{x}^{(N)}-\bar{x}) +Q(x_i^*-\bar{x}_i)-Q_{\Gamma} ({x}_*^{(N)}-\bar{x})\cr&-Q_{\Gamma_1} (x_0^*-\bar{x}_0)\big]+\tilde{x}_0^T\big[ F^T\Pi(\bar{x}^{(N)}-\bar{x})-Q_{\Gamma_1}^T(x^{(N)}_*-\bar{x})+\Gamma_1^TQ\Gamma_1(x_0^*-\bar{x}_0)\big]\Big\}dt\cr
&+ \frac{1}{N} \sum_{i=1}^N\mathbb{E}\big\{\tilde{x}_i^T(T)\big[H({x}^*_i(T)-\bar{x}_i(T))-H_{\bar{\Gamma}}({x}_*^{(N)}(T)-\bar{x}(T))-H_{\bar{\Gamma}_1} ({x}_0^*(T)-\bar{x}_0(T))\big]\cr&
 +\tilde{x}_0^T(T)\big[-H_{\bar{\Gamma}_1}^T({x}_*^{(N)}(T)-\bar{x}(T))+\bar{\Gamma}_1^TH\bar{\Gamma}_1(x_0^*(T)-\bar{x}_0(T))\big]\big\}.\cr
\leq&O(1/\sqrt{N})=\varepsilon_1.
\end{aligned}
\end{equation*}
Note that $\tilde{J}_i(u_0,\tilde{u})\geq0$.
Then we have
${J}_{\rm soc}(u_0,{u}^*)\leq J_{\rm soc}(u_0,u)+\varepsilon_1.$

 \emph{(For the leader).}  From (\ref{eq2}), we have
\begin{align}\label{eq36-a}
   J_0({u}_0^*,{u}^*)=&\mathbb{E}\int_{0}^{T}\big[\|\bar{x}_0^*- \Gamma_0\bar{x}+{x}_0^*-\bar{x}_0^*+\Gamma_0({x}^{(N)}_*-\bar{x})\|_{Q_0}^{2}+\|{u}_0^*\|_{R_0}^{2}\big]dt\cr
   &+\mathbb{E}\big[\|\bar{x}_0^*(T)- \bar{\Gamma}_0\bar{x}(T)+{x}_0^*(T)-\bar{x}_0^*(T)+\bar{\Gamma}_0({x}^{(N)}_*(T)-\bar{x}(T))\|_{H_0}^{2}\big]\cr
\leq& \bar{J}_0({u}_0^*,{u}^*)+\mathbb{E}\int_{0}^{T}\big[2\|{x}_0^*-\bar{x}_0^*\|^2_{Q_0}+2\|\Gamma_0({x}_*^{(N)}-\bar{x})\|_{Q_0}^{2}\big]dt\cr
   &+\int_{0}^{T}2\big(\mathbb{E}\|x_0^*- \Gamma_0\bar{x}\|^2\cdot2\|Q_0\Gamma_0\|^2(\mathbb{E}\|{x}_*^{(N)}-\bar{x}\|^2+\mathbb{E}\|{x}_*^{(N)}-\bar{x}\|^2)\big)^{1/2}dt\cr
   &+\mathbb{E}\big[2\|{x}_0^*(T)-\bar{x}_0^*(T)\|^2_{H_0}+2\|\bar{\Gamma}_0({x}_*^{(N)}(T)-\bar{x}(T))\|_{H_0}^{2}\big]\cr
   &+C\big[\mathbb{E}\|x_0^*(T)- \bar{\Gamma}_0\bar{x}(T)\|^2\cdot(\mathbb{E}\|{x}_*^{(N)}(T)-\bar{x}(T)\|^2+\mathbb{E}\|{x}_*^{(N)}(T)-\bar{x}(T)\|^2)\big]^{1/2}\cr
\leq &\bar{J}_0({u}_0^*,{u}^*)+O(1/\sqrt{N}).
\end{align}
  From  this and Theorem \ref{thm3.3}, we obtain
 \begin{equation}\label{eq41c}
   J_0({u}^*_0,{u}^*)\leq  \bar{J}_0(u_0,u^*)+O(1/\sqrt{N}).
    \end{equation}
From (\ref{eq39}),
\begin{align*}
   \bar{J}_0({u}_0,{u}^*)=&\mathbb{E}\int_{0}^{T}\big[\|\bar{x}_0
   -\Gamma_0 {x}^{(N)}_*+x_0-\bar{x}_0+\Gamma_0({x}_*^{(N)}-\bar{x})\|_{Q_0}^{2}+\|u_0\|_{R_0}^{2}\big]dt\cr
   &+\mathbb{E}\big[\|\bar{x}_0(T)
   - {x}^{(N)}_*(T)+x_0(T)-\bar{x}_0(T)+\bar{\Gamma}_0({x}_*^{(N)}(T)-\bar{x}(T))\|_{H_0}^{2}\big]\cr
\leq&{J}_0({u}_0,{u}^*)+\|\Gamma_0(x_0-\bar{x}_0+{x}_*^{(N)}-{x}^*)\|_{Q_0}^{2}\cr
   &+2\mathbb{E}\int_{0}^{T}\Big[\big(\|x_0- \Gamma_0{x}_*^{(N)}\|^2\|Q_0\Gamma_0(x_0-\bar{x}_0+{x}_*^{(N)}-\bar{x})\|^2\big)^{1/2}
   \Big]dt+O\big(\frac{1}{\sqrt{N}}\big)\cr
\leq &{J}_0({u}_0,{u}^*)+O(1/\sqrt{N}).
\end{align*}
From this and (\ref{eq41c}), we have
$J_0({u}^*_0,{u}^*)\leq {J}_0({u}_0,{u}^*)+\varepsilon_2,$
where $\varepsilon_2=O(1/\sqrt{N})$.
\hfill{$\Box$}

\section{Proof of Theorem \ref{thm3.5}.}\label{app3}
\def\theequation{C.\arabic{equation}}
\setcounter{equation}{0}

 To prove the theorem, we need a lemma.
Consider a MF-type system 
\begin{align*}
  d\bar{x}_i=& (A\bar{x}_i+B{u}_i+G\mathbb{E}_{\mathcal{F}^0}[\bar{x}_i]+F\bar{x}_0+B_1u_0^*)dt+DdW_i,\cr
  d\bar{x}_0=& (A_0\bar{x}_0+B_0{u}_0^*+G_0\mathbb{E}_{\mathcal{F}^0}[\bar{x}_i])dt+D_0dW_0
\end{align*}
with the cost function
\begin{align*}
  \mathcal{J}_i(u_i)=&\mathbb{E}\int_0^T\big(\|\bar{x}_i-\Gamma\mathbb{E}_{\mathcal{F}^0}[\bar{x}_i]-\Gamma_1\bar{x}_0\|^2_Q+\|u_i\|^2_R+2u_i^TLu_0^*
  +\|u_0^*\|^2_{R_1}\big)dt\cr
  &+\mathbb{E}\big[\|\bar{x}_i(T)-\Gamma\mathbb{E}_{\mathcal{F}^0}[\bar{x}_i(T)]-\Gamma_1\bar{x}_0(T)\|^2_H\big].
  \end{align*}

\begin{lemma}\label{lem3.2}
For the above problem, 
the optimal control is given by
$$u_i^*=-R^{-1}B^T[\Pi \bar{x}_i+(\bar{\Pi}-\Pi)\bar{x}+M\bar{x}_0+\bar{\varphi}]-R^{-1}Lu_0^*,$$ 
and the corresponding optimal cost is
$\mathbb{E}[\|\xi_i\|^2_{\Pi(0)}+\|\bar{\xi}\|^2_{\bar{\Pi}(0)-\Pi(0)}+\|\xi_0\|^2_{M_0(0)}+2\bar{\xi}^T\Pi_0(0)\xi_0]+m_T.$ 

\end{lemma}

\emph{Proof.} Note that $\mathbb{E}_{\mathcal{F}^0}[\bar{x}_i]=\bar{x}$ satisfies
\begin{align}\label{eq40d}
  d\bar{x}=& [(A+G)\bar{x}+B\bar{u}+F\bar{x}_0+B_1u_0^*]dt,
\end{align}
where $\bar{u}\stackrel{\Delta}{=}\mathbb{E}_{\mathcal{F}^0}[\bar{u}_i]$. Then we have
$$ d(\bar{x}_i-\bar{x})= [A(\bar{x}_i-\bar{x})+B({u}_i-\bar{u})]dt+DdW_i.$$
Applying It\^{o}'s formula to $\|\bar{x}_i-\bar{x}\|^2_{\Pi}$, it follows that
\begin{align}\label{eq55}
&\mathbb{E}\big[\|\bar{x}_i(T)-\bar{x}(T)\|^2_{H}-\|\bar{x}_i(0)-\bar{x}(0)\|^2_{{\Pi}{(0)}}\big]\cr
=&\mathbb{E}\int_0^T\big[(\bar{x}_i-\bar{x})^T(\dot{{\Pi}}+{\Pi}A+A^T{\Pi})(\bar{x}_i-\bar{x})+2(\bar{x}_i-\bar{x})^T{\Pi}B(u_i-\bar{u})+D^T\Pi D\big]dt.
\end{align}
Furthermore, 
by (\ref{eq40d}) and (\ref{eq4a}), one can obtain
\begin{align}
&\mathbb{E}[\|\bar{x}(T)\|^2_{H-H_{\bar{\Gamma}}}-\|\bar{x}(0)\|^2_{\bar{\Pi}{(0)}}]\cr
=&\mathbb{E}\int_0^T\big[\bar{x}^T(\dot{\bar{\Pi}}+\bar{\Pi}(A+G)+(A+G)^T\bar{\Pi})\bar{x}+2\bar{x}^T\bar{\Pi}(B\bar{u}+F\bar{x}_0+B_1u_0^*)\big]dt,\\
&\mathbb{E}[\|\bar{x}_0(T)\|^2_{\bar{\Gamma}_1^TH\bar{\Gamma}_1}-\|\bar{x}_0(0)\|^2_{M_0(0)}]\cr
=&\mathbb{E}\int_0^T\big[\bar{x}^T_0(\dot{M}_0+M_0A_0+A_0^TM_0)\bar{x}_0+2\bar{x}_0^TM_0(B_0{u}^*_0+G_0\bar{x})+D_0^TM_0D_0\big]dt
\end{align}
and
\begin{align}
&\mathbb{E}[-\bar{x}_0^T(T)H_{\bar{\Gamma}_1}^T\bar{x}(T)-\bar{x}_0^T(0){\Pi_0}(0)\bar{x}(0)]\cr
=&\mathbb{E}\int_0^T\big\{\bar{x}^T_0[\dot{\Pi}_0+A^T\Pi_0+\Pi_0(A+G)]\bar{x}\cr
&+\bar{x}^TG_0^T\Pi_0\bar{x}+({u}_0^*)^TB_0^T\Pi_0\bar{x}+\bar{x}^T_0\Pi_0B\bar{u}+\bar{x}_0^T\Pi_0F\bar{x}_0+\bar{x}_0^T\Pi_0B_1u_0^*
\big\}.
\end{align}
Also, applying
It\^{o}'s formula to $\bar{x}^T\bar{\varphi}$ and $\bar{x}^T_0\bar{\varphi}_0$, we have
\begin{align}
&\mathbb{E}[\bar{x}^T(T)\bar{\varphi}(T)-\bar{x}^T(0)\bar{\varphi}(0)]\cr
=&\mathbb{E}\int_0^T\big[\bar{x}^T\big(\bar{\Pi}BR^{-1}B^T\bar{\varphi}-G_0^T\bar{\varphi}_0-(\bar{\Pi}\bar{B}_1+MB_0){u}_0^*\big)+(B\bar{u}+F\bar{x}_0)^T\bar{\varphi}\big]dt,
\end{align}
and
\begin{align}\label{eq60}
\mathbb{E}[\bar{x}^T_0(T)\bar{\varphi}_0(T)-\bar{x}^T_0(0)\bar{\varphi}_0(0)]
=&\mathbb{E}\int_0^T\big[-\bar{x}_0^T(F^T-\bar{\Pi}_0BR^{-1}B^T\bar{\varphi}+({\Pi}_0\bar{B}_1+M_0B_0){u}_0^*)
\cr&+(B_0{u}_0^*+G_0\bar{x}_0)^T\bar{\varphi}_0+D_0^T(\bar{q}_0^0-M_0D_0)\big]dt.
\end{align}
Note that $\bar{x}=\mathbb{E}[\bar{x}_i|\mathcal{F}^0]$. From (\ref{eq55})-(\ref{eq60}), we obtain
\begin{align*}
  \mathcal{J}_i(u_i)=&\mathbb{E}\int_0^T\big[\|\bar{x}_i-\bar{x}\|^2_Q
  +2(\bar{x}_i-\bar{x})^TQ((I-\Gamma)\bar{x}-\Gamma_1\bar{x}_0)+\|(I-\Gamma)\bar{x}-\Gamma_1\bar{x}_0\|^2_Q\cr
  &+2u_i^TLu_0^*+\|u_i-\bar{u}\|^2_R+\|\bar{u}\|^2_R+\|u_0^*\|^2_{R_1}\big]dt+\mathbb{E}\big[\|\bar{x}_i(T)-\bar{\Gamma}\mathbb{E}[\bar{x}(T)]-\bar{\Gamma}_1\bar{x}_0(T)\|^2_H\big]\cr
  =&\mathbb{E}\int_0^T\big(\|\bar{x}_i-\bar{x}\|^2_Q+\|(I-\Gamma)\bar{x}-\Gamma_1\bar{x}_0\|^2_Q
  +\|u_i-\bar{u}\|^2_R+2\bar{u}^TLu_0^*+\|\bar{u}\|^2_R+\|u_0^*\|^2_{R_1}\big)dt\cr
 & +\mathbb{E}\big[\|\bar{x}_i(T)-\bar{x}(T)\|^2_H+\|(I-\bar{\Gamma})\bar{x}(T)-\bar{\Gamma}_1\bar{x}_0(T)\|^2_H\big]
 \end{align*}
 \begin{align*}
   =& \mathbb{E}\big[\|\bar{x}_i(0)-\bar{x}(0)\|^2_{\Pi(0)}+\|\bar{x}(0)\|^2_{\bar{\Pi}(0)}+\|\bar{x}_0(0)\|^2_{M_0(0)}+2\bar{x}_0^T(0)\Pi_0(0)\bar{x}(0)+2\bar{x}^T(0)\bar{\varphi}(0)
  \cr
  &+2\bar{x}^T_0(0)\bar{\varphi}_0(0)  \big]+\mathbb{E}\int_0^T\big[(\bar{x}_i-\bar{x})^T\Pi BR^{-1}B^T\Pi(\bar{x}_i-\bar{x}) +\bar{x}^T\bar{\Pi} BR^{-1}B^T\bar{\Pi}\bar{x}\cr
  &+\bar{x}^T_0M^T BR^{-1}B^TM \bar{x}_0
  +\bar{x}^T_0M^T BR^{-1}B^T\bar{\Pi} \bar{x}+2(x_i-\bar{x})^T{\Pi}B(u_i-\bar{u})
  \cr
  &+2\bar{x}^T\bar{\Pi}B\bar{u}+2\bar{x}^T_0\Pi_0B\bar{u}+2\bar{\varphi}^TB\bar{u}+2\bar{x}^T\bar{\Pi}BR^{-1}B^T\bar{\varphi}
+2\bar{x}_0^TM^TBR^{-1}B^T\bar{\varphi} \cr  &
+\|u_i-\bar{u}\|^2_R+\|\bar{u}\|^2_R+2\bar{u}^TLu_0^*+\|u_0^*\|^2_{R_1}+2(B^T\bar{\Pi}\bar{x}+B^T{\Pi}_0^T\bar{x}_0)^TR^{-1}Lu_0^*\cr
&+2\bar{\varphi}_0^TB_0{u}_0^*+D^T\Pi D+D_0^TM_0D_0+D_0^T(\bar{q}_0^0-M_0D_0)\big]dt\cr
=& \mathbb{E}[\|\xi_i-\bar{\xi}\|^2_{\Pi(0)}+\|\bar{\xi}\|^2_{\bar{\Pi}(0)}+\|\xi_0\|^2_{M_0(0)}+2\bar{\xi}^T\Pi_0(0)\xi_0 \big]+m_T\cr&
+\mathbb{E}\int_0^T\big[\|u_i-\bar{u}+R^{-1}B^T(\bar{x}_i-\bar{x})\|^2_R+\|\bar{u}+R^{-1}B^T(\bar{\Pi}\bar{x}+M\bar{x}_0+\bar{\varphi})+R^{-1}Lu_0^*\|^2_R\big]dt\cr
\geq& \mathbb{E}[\|\xi_i\|^2_{\Pi(0)}+\|\bar{\xi}\|^2_{\bar{\Pi}(0)-\Pi(0)}+\|\xi_0\|^2_{M_0(0)}+2\bar{\xi}^T\Pi_0(0)\xi_0 \big]+m_T. 
  \end{align*}
  \rightline{$\Box$}
  \emph{Proof of Theorem \ref{thm3.5}.} (For followers) 
   Note that $u_i^*=\!-R^{-1}B^T(\Pi \bar{x}_i+(\bar{\Pi}-\Pi)\bar{x}+M\bar{x}_0+\bar{\varphi})\!-\!R^{-1}Lu_0^*$, and
   \begin{align*}
        \frac{1}{N}J_{\rm soc}^{(N)}(u^*,u_0^*)
        =&\frac{1}{N}\sum_{i=1}^N\mathbb{E}\int_0^T\big[\|\bar{x}_i-\Gamma \bar{x}-\Gamma_1\bar{x}_0+x_i^*-\bar{x}_i-\Gamma (x_*^{(N)}-\bar{x})-\Gamma_1(x_0^*-\bar{x}_0)\|^2_Q\cr
        & +\|u_i^*\|^2_{R}+2u_i^*Lu_0^*+\|u_0^*\|^2_{R_1}\big]dt+\mathbb{E}[\|x_i^*(T)-\bar{\Gamma} x_*^{(N)}(T)-\bar{\Gamma}_1x_0^*(T)\|^2_H].
  \end{align*}
 By Schwarz's inequality and Lemma \ref{lem3.1}, one can obtain
    \begin{align*}
  &\frac{1}{N}\big|J_{\rm soc}^{(N)}(u^*,u_0^*)-\sum_{i=1}^N \mathcal{J}_i(u_i^*)  \big|\cr
  \leq& \frac{1}{N}\sum_{i=1}^N\mathbb{E}\int_0^T\big[\|x_i^*-\bar{x}_i\|^2_Q+\|\Gamma (x_*^{(N)}-\bar{x})\|^2_Q+\|\Gamma_1(x_0^*-\bar{x}_0)\|^2_Q
  \big]dt\cr
  &+\frac{C}{N}\sum_{i=1}^N\sup_{0\leq t\leq T}\big(\mathbb{E}\|x_i^*-\bar{x}_i\|^2_Q\big)^{1/2}+\frac{C}{N}\sum_{i=1}^N\sup_{0\leq t\leq T}\big(\mathbb{E}\|\Gamma (x_*^{(N)}-\bar{x})\|^2_Q\big)^{1/2}\cr
  &+\frac{C}{N}\sum_{i=1}^N\sup_{0\leq t\leq T}\big(\mathbb{E}\|\Gamma_1(x_0^*-\bar{x}_0)\|^2_Q\big)^{1/2}\leq O(1/\sqrt{N}).
    \end{align*}
    From this and Lemma \ref{lem3.2}, we have (\ref{eq46a}).

    (For the leader) By a similar argument with the proof of Theorem \ref{thm3.3}, one can obtain
        $$\begin{aligned}
           \bar{J}_{0}(u_0^*,u^*)=&\mathbb{E}\Big\{{\xi}^T_0y_0(0)+\bar{\xi}^T\bar{y}(0)+\int_0^T\big[\big\langle R_0u_0^*+B_0^Ty_0+\bar{B}_1^T\bar{y}-({\Pi}_0\bar{B}_1+M_0B_0)^T\psi_0\cr
  &-(\bar{\Pi}\bar{B}_1+MB_0)^T\psi, u_0^*\big\rangle +D_0^T\beta_0\big]dt\Big\}\cr
  =&\mathbb{E}\Big[{\xi}^T_0y_0(0)+\bar{\xi}^T\bar{y}(0)+\int_0^T\big(D_0^T\beta_0\big)dt\Big].
          \end{aligned}$$
          By (\ref{eq36-a}), we have (\ref{eq46}). The proof of the theorem is completed.
    \hfill{$\Box$}

    \section{Proofs of Theorems \ref{thm4.1} and \ref{thm4.3}.}\label{app4}
\def\theequation{D.\arabic{equation}}
\setcounter{equation}{0}

{\it Proof.}
Suppose that $\breve{u}_i=-R^{-1}B^T{p}_i-R^{-1}L(P_0x_0+\bar{P}x^{(N)}),$ where {${({p}_i, {q}_{i}^{j}},i,j=0,1,\cdots,N)$ is a set of solutions to in (\ref{eq3d}).}
 Denote by $\breve{x}_0$, $\breve{x}_i$ the state of agent $i$ under the control $\breve{u}_i$, $i=1,\cdots,N$. For any $u_i\in L^2_{{\mathcal  F}}(0, T; \mathbb{R}^r) $
 and $\lambda\in \mathbb{R}\ (\lambda\not=0)$, let $u_i^{\lambda}=\breve{u}_i+\lambda u_i$.
 Denote by $ x_i^{\lambda}$, $i=0,1,\cdots,N$ the solution of the following perturbed state equation
$$\left\{ \begin{aligned}
 dx_0^{\lambda}(t)=&\Big[(A_0+B_0P_0)x_0^{\lambda}(t)
+\frac{1}{N}(G_0+B_0\bar{P})\sum_{i=1}^Nx^{\lambda}_i(t)\Big]dt
+D_0dW_0(t),
\\
dx_i^{\lambda}(t)=&\Big[Ax_i^{\lambda}(t)+B(\breve{u}_i(t)+\lambda u_i(t))+\frac{1}{N}(G+B_1\bar{P})\sum_{i=1}^Nx^{\lambda}_i(t)\cr
&+(F+B_1P_0)x^{\lambda}_0(t)\Big]dt+DdW_i(t),\cr
x_0(0)=&\xi_{0},\ x_i(0)=\xi_{i},\ i=1,2,\cdots,N.
 \end{aligned}\right.$$
 Let $z_i=(x_i^{\lambda}-\breve{x}_i)/\lambda$, $i=0,1,2,\cdots,N$. 
It can be verified that
 $z_0$ and $z_i$ satisfy
 $$\left\{ \begin{aligned}
 dz_0(t)=&\big[(A_0+B_0P_0)z_0(t)+(G_0+B_0\bar{P})z^{(N)}(t)\big]dt,
\\
dz_i(t)=&\big[Az_i(t)+B u_i(t)+(G+B_1\bar{P})z^{(N)}(t)+(F+B_1P_0)z_0(t)\big]dt,\cr
z_0(0)=&\ 0,\ z_i(0)=0,\ i=1,2,\cdots,N,
 \end{aligned}\right.$$
where $z^{(N)}=\frac{1}{N}\sum_{i=1}^Nz_i$. Then by It\^{o}'s formula, we have
{\begin{align}\label{eq6-a}
&\mathbb{E}[\langle- H_{\bar{\Gamma}_1} ^Tx^{(N)}(T)+\bar{\Gamma}_1^TH\bar{\Gamma}_1x_0(T),z_0(T)\rangle ]=\mathbb{E}[\langle \breve{p}_0(T),z_0(T)\rangle-\langle \breve{p}_0(0),z_0(0)\rangle]\cr
=&\mathbb{E}\!\int_0^T\! \Big\{\big\langle -\big[(A_0+B_0P_0)^T\breve{p}_0(t)+(F+B_1P_0)^T\breve{p}^{(N)}(t)+P_0^TR_1(P_0\breve{x}_0(t)+\bar{P}\breve{x}^{(N)}(t))\cr
&-\Gamma_1Q\big((I-\Gamma) \breve{x}^{(N)}(t)-\Gamma_1 \breve{x}_0(t)\big)+P_0^TL^T\breve{u}^{(N)}(t)\big],
z_0(t)\rangle\cr
&+\langle \breve{p}_0(t),(A_0+B_0P_0)z_0(t)+(G_0+B_0\bar{P})z^{(N)}(t)\rangle\Big\} dt\cr
=&\mathbb{E}\!\int_0^T\!  \Big\{\big\langle -\big[(F+B_1P_0)^T\breve{p}^{(N)}(t)+P_0^TR_1(P_0\breve{x}_0(t)+\bar{P}\breve{x}^{(N)}(t))+P_0^TL^T\breve{u}^{(N)}(t)\cr
&-\Gamma_1Q\big((I-\Gamma) \breve{x}^{(N)}(t)-\Gamma_1 \breve{x}_0(t)\big)\big],
z_0(t)\rangle+\langle \breve{p}_0(t),(G_0+B_0\bar{P})z^{(N)}(t)\rangle \Big\} dt,
\end{align}}
and
{\begin{align}\label{eq7-b}
&\sum_{i=1}^N\mathbb{E}[\langle Hx_i(T)-H_{\bar{\Gamma}}x^{(N)}(T)-H_{\bar{\Gamma}_1} x_0(T),z_i(T)\rangle]
=\sum_{i=1}^N\mathbb{E}[\langle \breve{p}_i(T),z_i(T)\rangle-\langle \breve{p}_i(0),z_i(0)\rangle]\cr
=&\sum_{i=1}^N\mathbb{E}\!\int_0^T\!  \Big\{\big\langle -\big[A^T\breve{p}_i(t)+(G+B_1\bar{P})^T\breve{p}^{(N)}(t)+(G_0+B_0\bar{P})^T\breve{p}_0(t)+Q\breve{x}_i(t)
\cr&+(\bar{P}^TR_1\bar{P}-Q_{\Gamma})\breve{x}^{(N)}(t)
+\big((\Gamma-I)^TQ\Gamma_1+\bar{P}^TR_1P_0\big)\breve{x}_0(t)+\bar{P}^TL^T\breve{u}^{(N)}(t)\big],z_i(t)\big\rangle\cr
&+\langle \breve{p}_i(t),Az_i(t)+Bu_i(t)+(G+B_1\bar{P})z^{(N)}(t)+(F+B_1P_0)z_0(t)\rangle \Big\} dt\cr
= &\sum_{i=1}^N \mathbb{E}\int_0^T \Big\{\big\langle -\big[(G_0+B_0\bar{P})^T\breve{p}_0(t)+Q\breve{x}_i(t)
+(\bar{P}^TR_1\bar{P}-Q_{\Gamma})\breve{x}^{(N)}(t)
\cr&+\big((\Gamma-I)^TQ\Gamma_1+\bar{P}^TR_1P_0\big)\breve{x}_0(t)+\bar{P}^TL^T\breve{u}_i(t)\big],z_i(t)\big\rangle\cr
&+\langle (F+B_1P_0)^T\breve{p}^{(N)}(t),z_0(t)\rangle+\langle B^T\breve{p}_i(t),\breve{u}^{(N)}(t)\rangle\Big\} dt.
\end{align}}
From (\ref{eq3}),
$$
\begin{aligned}
 &{J}_{\rm soc}(\breve{u}+\lambda u)-{J}_{\rm soc}(\breve{u})=2\lambda I_1+{\lambda^2}I_2
\end{aligned}
$$
where $\breve{u}=(\breve{u}_1,\cdots,\breve{u}_N)$, and  (suppressing the time $t$)
\begin{align*}
I_1\stackrel{\Delta}{=}&\sum_{i=1}^N\mathbb{E}\int_0^T \big[\big\langle Q\big(\breve{x}_i-(\Gamma\breve{x}^{(N)}+\Gamma_1\breve{x}_0)\big), z_i-\Gamma z^{(N)}-\Gamma_1z_0\big\rangle
+\langle R \breve{u}_i,u_i\rangle \cr&+\langle L(P_0\breve{x}_0+\bar{P}\breve{x}^{(N)}), u_i\rangle
  +\langle L(P_0z_0+\bar{P}z^{(N)}),\breve{u}_i\rangle
+\langle R_1(P_0\breve{x}_0+\bar{P}\breve{x}^{(N)}),P_0z_0+\bar{P}z^{(N)}\rangle\big]dt\cr &
+\mathbb{E} \big[\big\langle {H}(\breve{x}_i(T)
   -\bar{\Gamma} \breve{x}^{(N)}(T)-\bar{\Gamma}_1\breve{x}_0(T)),z_i(T)
   -\bar{\Gamma} z^{(N)}(T)-\bar{\Gamma}_1z_0(T)\big\rangle\big],
\cr
I_2\stackrel{\Delta}{=}&\sum_{i=1}^N\mathbb{E}\int_0^T\! \big[\big\|z_i
  \! -\!\Gamma z^{(N)}-\Gamma_1z_0\big\|^2_{Q}
+\|u_i\|^2_{R}+2u_i^TL(P_0z_0+\bar{P}z^{(N)})\cr
&+\|P_0z_0+\bar{P}z^{(N)}\|^2_{R_1}\big]dt+\mathbb{E} \big[\big\|z_i(T)
   -\bar{\Gamma} z^{(N)}(T)-\bar{\Gamma}_1z_0(T)\big\|^2_{H}\big].
\end{align*}
Note that
\begin{align*}
  &\sum_{i=1}^N\mathbb{E}\int_0^T \big\langle Q\big(\breve{x}_i-(\Gamma\breve{x}^{(N)}+\Gamma_1\breve{x}_0)\big),\Gamma z^{(N)}\big\rangle dt\cr
  &+\sum_{i=1}^N\mathbb{E} \big[\big\langle {H}(\breve{x}_i(T)
   -\bar{\Gamma} \breve{x}^{(N)}(T)-\bar{\Gamma}_1\breve{x}_0(T)),\bar{\Gamma} z^{(N)}(T)\big\rangle\big]\cr
 =& \sum_{j=1}^N \mathbb{E}\int_0^T \big\langle  {\Gamma^TQ} \big((I-\Gamma)\breve{x}^{(N)}-\Gamma_1\breve{x}_0\big), z_j\big\rangle  dt\cr
 &
 +\sum_{j=1}^N\mathbb{E} \big[\big\langle\bar{\Gamma}^T {H}((I-
   \bar{\Gamma}) \breve{x}^{(N)}(T)-\bar{\Gamma}_1\breve{x}_0(T)),z_j(T)\big\rangle\big],
\end{align*}
and
$$\sum_{i=1}^N \mathbb{E}\int_0^T \langle L\bar{P}z^{(N)},\breve{u}_i\rangle dt=\sum_{i=1}^N \mathbb{E}\int_0^T \langle\bar{P}^TL^T\breve{u}^{(N)},z_i\rangle dt.$$
From (\ref{eq6-a}) and (\ref{eq7-b}), one can obtain that
{\begin{align}\label{eq10b}
I_1=&\sum_{i=1}^N\mathbb{E}\int_0^T\big[\big\langle Q\big(\breve{x}_i-(\Gamma\breve{x}^{(N)}+\Gamma_1\breve{x}_0)\big), z_i-\Gamma z^{(N)}-\Gamma_1z_0\big\rangle\cr
&+\langle R \breve{u}_i,u_i\rangle+\langle L(P_0\breve{x}_0+\bar{P}\breve{x}^{(N)}), u_i\rangle
  +\langle L(P_0z_0+\bar{P}z^{(N)}),\breve{u}_i\rangle
\cr &+\langle R_1(P_0\breve{x}_0+\bar{P}\breve{x}^{(N)}),P_0z_0+\bar{P}z^{(N)}\rangle\big]dt\cr
&+\sum_{i=1}^N \mathbb{E}\int_0^T  \big\{\langle-\big[P_0^TR_1(P_0\breve{x}_0+\bar{P}\breve{x}^{(N)})+P_0^TL^T\breve{u}^{(N)}\cr
&-\Gamma_1Q\big((I-\Gamma) \breve{x}^{(N)}-\Gamma_1 \breve{x}_0\big)\big],
z_0\rangle+\langle \breve{p}_0,(G_0+B_0\bar{P})z^{(N)}\rangle\cr
& -\big[(G_0+B_0\bar{P})^T\breve{p}_0+Q\breve{x}_i
+(\bar{P}^TR_1\bar{P}-Q_{\Gamma})\breve{x}^{(N)}
\cr&+\big((\Gamma-I)^TQ\Gamma_1+\bar{P}^TR_1P_0\big)\breve{x}_0+\bar{P}^TL^T\breve{u}^{(N)}\big],z_i\big\rangle+\langle B^T\breve{p}_i,u_i\rangle \big\}dt\cr
=&\sum_{i=1}^N\mathbb{E}\int_0^T\big[\big\langle Q\breve{x}_i-Q_{\Gamma}\breve{x}^{(N)}+(\Gamma-I)^TQ\Gamma_1\breve{x}_0)+\bar{P}^TR_1(P_0\breve{x}_0+\bar{P}\breve{x}^{(N)}), z_i\big\rangle\cr
&-\big\langle \Gamma_1^TQ\big((I-\Gamma)\breve{x}^{(N)}-\Gamma_1\breve{x}_0\big)+P_0^TR_1(P_0\breve{x}_0+\bar{P}\breve{x}^{(N)}),z_0\big\rangle \cr
&+
\langle R\breve{u}_i+ L(P_0\breve{x}_0+\bar{P}\breve{x}^{(N)})+ B^T\breve{p}_i,u_i\rangle \big]dt\cr
&+\sum_{i=1}^N \mathbb{E}\int_0^T  \big\{\langle-\big[P_0^TR_1(P_0\breve{x}_0+\bar{P}\breve{x}^{(N)})-\Gamma_1Q\big((I-\Gamma) \breve{x}^{(N)}-\Gamma_1 \breve{x}_0\big)\big],
z_0\rangle\cr
& -\big[Q\breve{x}_i
+(\bar{P}^TR_1\bar{P}-Q_{\Gamma})\breve{x}^{(N)}
+\big((\Gamma-I)^TQ\Gamma_1+\bar{P}^TR_1P_0\big)\breve{x}_0\big],z_i\big\rangle \big\}dt\cr
=&\sum_{i=1}^N\mathbb{E}\int_0^T\big\langle R\breve{u}_i+ L(P_0\breve{x}_0+\bar{P}\breve{x}^{(N)})+B^T\breve{p}_i,u_i\big\rangle dt.
\end{align}}
Note that
\begin{align*}
I_2{=}&\sum_{i=1}^N\mathbb{E}\int_0^T\! \big[\big\|z_i
  \! -\!\Gamma z^{(N)}-\Gamma_1z_0\big\|^2_{Q}
+\langle (R-LR_1^{-1}L^T)u_i,u_i\rangle\cr
&+\|R_1^{-1}L^Tu_i+P_0z_0+\bar{P}z^{(N)}\|^2_{R_1}\big]dt+\mathbb{E}\big[\big\|z_i(T)
   -\bar{\Gamma} z^{(N)}(T)-\bar{\Gamma}_1z_0(T)\big\|^2_{H}\big].
\end{align*}
Since $Q\geq0$, $R>0$ and $R-LR_1^{-1}L^T\geq 0$, we have $I_2\geq0$ and Problem (P3) admits an optimal control \cite{YZ99}. From (\ref{eq5a}), $\breve{u}$ is a minimizer to Problem (P1) if and only if
 $I_1=0 $, which is equivalent to
$\breve{u}_i=-{R^{-1}}\big(B^T\breve{p}_i+ L(P_0\breve{x}_0+\bar{P}\breve{x}^{(N)})\big).$
Thus, we have the optimality system (\ref{eq3d}).
This implies that (\ref{eq3d}) admits a solution $(\breve{x}_i,\breve{p}_i,\breve{q}_{i}^{j}, i,j=1,\cdots,N)$. \hfill{$\Box$}

\emph{Proof of Theorem \ref{thm4.3}.} \emph{(For followers).} Suppose $u_0=P_0x_0+\bar{P}\bar{x}$.  From (\ref{eq1}),
we have
\begin{equation}\label{eq71}
  dx^{(N)}=[(A+G)x^{(N)}+Bu^{(N)}+Fx_0+B_1(P_0x_0+\bar{P}\bar{x})]dt+\frac{1}{N}\sum_{j=1}^NDdW_j,
\end{equation}
where $u^{(N)}=\frac{1}{N}\sum_{j=1}^N u_j$. This leads to
\begin{equation*}
  d(x_i-x^{(N)})=[A(x_i-x^{(N)})+B(u_i-u^{(N)})]dt+DdW_i-\frac{1}{N}\sum_{j=1}^NDdW_j.
\end{equation*}
Applying It\^{o}'s formula to $\|x_i-x^{(N)}\|^2_{{K}}$,
\begin{align}\label{eq72}
  &\mathbb{E}[\|x_i(T)-x^{(N)}(T)\|^2_{H}-\|x_i(0)-x^{(N)}(0)\|^2_{K(0)}]\cr
  =&\mathbb{E}\int_0^T\Big\{(x_i-x^{(N)})^T (\dot{K}+
A^T{K}+{K}A  )(x_i-x^{(N)})\cr
&+2(x_i-x^{(N)})^TKB(u_i-u^{(N)})+\frac{N-1}{N}D^T{K}D\Big\}
dt.
\end{align}
Noting that $u_0=P_0x_0+\bar{P}\bar{x}$, 
 we have
$$dx_0=[(A_0+B_0P_0)x_0+G_0x^{(N)}+B_0\bar{P}\bar{x}]dt+D_0dW_0.$$
From 
this and (\ref{eq71}), 
 we obtain
\begin{align}\label{eq73}
  &\mathbb{E}\big[\|x^{(N)}(T)\|^2_{H-H_{\bar{\Gamma}}}-\|x^{(N)}(0)\|^2_{K(0)+\bar{K}(0)}\big]\cr
  =&\mathbb{E}\int_0^T\Big\{(x^{(N)})^T[\dot{K}+\dot{\bar{K}}+
(A+G)^T(K+\bar{K})+(K+\bar{K})(A+G)  ]x^{(N)}+2(x^{(N)})^T(K+\bar{K})u^{(N)}\cr
&+2(x^{(N)})^T(K+\bar{K})[(F+B_1P_0)x_0+B_1\bar{P}\bar{x}]
+\frac{1}{N}D^T(K+\bar{K})D\Big\}
dt
\end{align}
and
\begin{align}\label{eq74}
  \mathbb{E}[\|x_0(T)\|^2_{\bar{\Gamma}_1^TH\bar{\Gamma}_1}-\|x_0(0)\|^2_{\Lambda_0(0)}]
  =&\mathbb{E}\int_0^T\Big\{x_0^T[\dot{\Lambda}_0+
(A_0+B_0P_0)^T\Lambda_0+\Lambda_0(A_0+B_0P_0) ]x_0\cr
&+2(G_0 x^{(N)}+B_0\bar{P}\bar{x})^T\Lambda_0x_0+D_0^T\Lambda_0D_0\Big\}
dt.
\end{align}
Applying It\^{o}'s formula to $x_0^T{\bar{\Lambda}}x^{(N)}$ and $(x^{(N)})^TK_0x_0$, we have
\begin{align}\label{eq75}
  &\mathbb{E}[x_0^T(T)({-H_{\bar{\Gamma}_1}}^T)x^{(N)}(T)-x_0^T(0){\bar{\Lambda}}(0)x^{(N)}(0)]\cr
  =&\mathbb{E}\int_0^T\Big\{x_0^T[\dot{\bar{\Lambda}}+\bar{\Lambda}
(A+G)+(A_0+B_0P_0)^T\bar{\Lambda} ]x^{(N)}+x_0^T\bar{\Lambda}Bu^{(N)}\cr
&+x_0^T\bar{\Lambda}(F+B_1P_0)x_0+x_0^T\bar{\Lambda}B_1\bar{P}\bar{x}+(G_0x^{(N)}+B_0\bar{P}\bar{x})^T\bar{\Lambda}x^{(N)}\Big\}
dt
\end{align}
and
\begin{align}\label{eq76}
  &\mathbb{E}[(x^{(N)}(T))^T(-H_{\bar{\Gamma}_1})x_0(T)-(x^{(N)}(0))^T{K_0}(0)x_0(0)]\cr
  =&\mathbb{E}\int_0^T\Big\{(x^{(N)})^T[\dot{K}_0+K_0(A_0+B_0P_0)+
(A+G)^TK_0 ]x_0+(u^{(N)})^TB^TK_0 x_0 \cr
&+(x^{(N)})^TK_0(G_0x^{(N)}+B_0\bar{P}\bar{x})+x_0^T(F+B_1P_0)^TK_0x_0+\bar{x}^T \bar{P}^TB_1^TK_0x_0\Big\}
dt.
\end{align}
By (\ref{eq3c}), we have
$$
d(x^{(N)}-\bar{x})=(A+G+BK)(x^{(N)}-\bar{x})dt+\frac{1}{N}\sum_{i=1}^NDdW_i,
$$
which implies
\begin{equation}\label{eq91}
\sup_{0\leq t\leq T}\mathbb{E}\|x^{(N)}-\bar{x}\|^2=O(1/N).
\end{equation}
From (\ref{eq32a}) and (\ref{eq72})-(\ref{eq76}), for any $u=\{u_i\in \mathcal{U}_c,i=1,\cdots,N\}$,
\begin{align}
\label{eq87}
  &J_{\rm soc}(u,u_0)\cr
  =&\sum_{i=1}^N\mathbb{E}\int_0^T \big[\|x_i\|^2_Q-\|x^{(N)}\|^2_{Q_{\Gamma}}-2(x_i-\Gamma x^{(N)})^TQ\Gamma_1x_0+
  \|\Gamma_1x_0\|_Q^2\cr
 &+\|u_i\|^2_R+ 2u_i^TL(P_0x_0+\bar{P}\bar{x})+\|P_0x_0+\bar{P}\bar{x}\|^2_{R_1}\big]dt\cr
 &+\sum_{i=1}^N \mathbb{E}\big[\|x_i(T)\|_H^2-\|{x}^{(N)}(T)\|^2_{H_{\bar{\Gamma}}}
 -2(x_i(T)-\bar{\Gamma} x^{(N)}(T))^TH\bar{\Gamma}_1x_0(T) +\|\Gamma_1x_0(T)\|^2_H\big]\cr
  =&\sum_{i=1}^N\mathbb{E}\int_0^T \big[\|x_i-x^{(N)}\|^2_Q+\|x^{(N)}\|^2_{Q-Q_{\Gamma}}+2(x^{(N)})^T(\Gamma-I)Q\Gamma_1x_0+
  \|\Gamma_1x_0\|_Q^2\cr
 & +\|u_i-u^{(N)}\|^2_R+\|u^{(N)}\|^2_R+ 2(u^{(N)})^TL(P_0x_0+\bar{P}\bar{x})+\|P_0x_0+\bar{P}\bar{x}\|^2_{R_1}\big]dt\cr
  &+\sum_{i=1}^N \mathbb{E}\big[\|x_i(T)-{x}^{(N)}(T)\|_H^2+\|{x}^{(N)}(T)\|^2_{H-H_{\bar{\Gamma}}}
 -2[ x^{(N)}(T)]^TH_{\bar{\Gamma}_1}x_0(T) +\|\Gamma_1x_0(T)\|^2_H\big]
 \end{align}
 \begin{align*}
  =&\sum_{i=1}^N\mathbb{E}[\|x_i(0)-x^{(N)}(0)\|^2_{{K(0)}}+\|x^{(N)}(0)\|^2_{K(0)+\bar{K}(0)}+\|x_0(0)\|^2_{\Lambda_0(0)}+2x_0^T(0){\bar{\Lambda}}x^{(N)}(0)]\cr
  &+\mathbb{E}\int_0^T\Big\{(x_i-x^{(N)})^TKBR^{-1}B^TK(x_i-x^{(N)})+2(x_i-x^{(N)})^TKB(u_i-u^{(N)})\cr
  &  +(x^{(N)})^T(K+\bar{K})BR^{-1}B^T(K+\bar{K})x^{(N)}+2(x^{(N)})^T(K+\bar{K})Bu^{(N)}\cr
  &+x_0^TK_0^TBR^{-1}B^TK_0 x_0+2(x^{(N)})^T(K+\bar{K})BR^{-1}B^TK_0x_0+2x_0^TK_0^TBu^{(N)}\cr
  &+2(u^{(N)})^TL(P_0x_0+\bar{P}\bar{x})+\|P_0x_0+\bar{P}\bar{x}\|^2_{R_1}+\|P_0x_0+\bar{P}{x}^{(N)}\|^2_{L^TR^{-1}L-R_1}
\cr & +2(x^{(N)})^T\bar{P}^TLR^{-1}B^T(K+\bar{K})x^{(N)}+2(x^{(N)})^T\bar{P}^TL^TR^{-1}B^TK_0x_0\cr
&+2x_0^TP_0^TL^TR^{-1}B^T(K+\bar{K})x^{(N)}+2x_0^TP_0L^TR^{-1}B^TK_0x_0+D^T{K}D+\frac{1}{N}D^T\bar{K}D\cr
&+D_0^T\Lambda_0D_0-2(x^{(N)}-\bar{x})^T\bar{P}^T\big[B_0^T(\Lambda_0x_0+\bar{\Lambda}x^{(N)})
+B_1^T(K+\bar{K})x^{(N)}+B^T_1K_0x_0\big]\Big\}dt\cr
=&\sum_{i=1}^N\mathbb{E}\Big\{[\|\xi_i\|^2_{K(0)}+\|\xi^{(N)}\|^2_{\bar{K}(0)} +\|\xi_0\|^2_{\Lambda_0}+2\xi_0^T{\bar{\Lambda}(0)}\xi^{(N)}+\int_0^T\Big[\|u_i-u^{(N)}+R^{-1}B^TK(x_i-x^{(N)})\|_R^2 \cr
&+\|u^{(N)}+R^{-1}B^T(K_0x_0+(K+\bar{K})x^{(N)})+R^{-1}L(P_0x_0+\bar{P}x^{(N)})\|_R^2\cr
&+D^T{K}D+\frac{1}{N}D^T\bar{K}D+D_0^T\Lambda_0D_0-2(x^{(N)}-\bar{x})^T\bar{P}^T\big[L^Tu^{(N)}\cr
&+R_1\big(P_0x_0+\frac{1}{2}\bar{P}\bar{x}+\frac{1}{2}\bar{P}x^{(N)}\big)+B_0^T(\Lambda_0x_0+\bar{\Lambda}x^{(N)})
+B_1^T(K+\bar{K})x^{(N)}+B^T_1K_0x_0\big]\Big]dt\Big\}\cr
\geq &\sum_{i=1}^N\mathbb{E}[\|\xi_i\|^2_{K}]+N\big(\mathbb{E}[\|\xi^{(N)}\|^2_{\bar{K}(0)} +\|\xi_0\|^2_{\Lambda_0}+2\xi_0^T{\bar{\Lambda}(0)}\xi^{(N)}]+D^T{K}D+D_0^T\Lambda_0D_0\big)\cr
&+D^T\bar{K}D-2N\Big(\mathbb{E}\int_0^T\|x^{(N)}-\bar{x}\|^2dt\cdot\mathbb{E}\int_0^T\|\bar{P}\|^2\big\|L^Tu^{(N)}+R_1(2P_0x_0+\bar{P}\bar{x}+\bar{P}x^{(N)})\cr
&+B_0^T(\Lambda_0x_0+\bar{\Lambda}x^{(N)})+B_1^T(K+\bar{K})x^{(N)}+B^TK_0x_0\big\|^2dt\Big)^{1/2}.
\end{align*}
Particularly, if $u_i=\hat{u}_i=-R^{-1}B^T(Kx_i+\bar{K}\bar{x}+K_0x_0)-R^{-1}L(P_0x_0+\bar{P}x^{(N)})$, then by (\ref{eq91}),
 we obtain (\ref{eq77}).
This implies $\hat{u}$ in (\ref{eq33a}) has $\epsilon_1$-social optimality,
where $\epsilon_1=O(1/\sqrt{N})$.

  \emph{(For the leader).} 
From (\ref{eq2}) and (\ref{eq91}), we have
\begin{align}\label{eq36a}
   J_0(\hat{u}_0,\hat{u})=&\mathbb{E}\int_{0}^{T}\big[\|x_0- \Gamma_0\bar{x}-\Gamma_0(\hat{x}^{(N)}-\bar{x})\|_{Q_0}^{2}+\|\hat{u}_0\|_{R_0}^{2}\big]dt
   +\mathbb{E}\big[\|x_0(T)- \bar{\Gamma}_0\hat{x}^{(N)}(T)\|_{H_0}^{2}\big]\cr
\leq& \bar{J}_0(\hat{u}_0,\hat{u})+C\sup_{0\leq t\leq T}\Big[2\big(\mathbb{E}\|x_0- \Gamma_0\bar{x}\|^2\mathbb{E}\|(\hat{x}^{(N)}-\bar{x})\|^2\big)^{1/2}+\mathbb{E}\|\hat{x}^{(N)}-\bar{x}\|^{2}\Big]dt\cr
\leq &\bar{J}_0(\hat{u}_0,\hat{u})+O(1/\sqrt{N}).
\end{align}
By It\^{o}'s formula,
  \begin{align}\label{eq32b}
  &\mathbb{E}[\bar{x}_0^T(T)H_0\bar{x}_0(T)]-  \mathbb{E}[\bar{x}_0^T(0)\Psi_1(0)\bar{x}_0(0)]\cr
  =&\mathbb{E} \int_0^T\big[\bar{x}_0^T(\dot{\Psi}_1+A_0^T\Psi_1+\Psi_1A_0)\bar{x}_0+2\bar{x}_0^T\Psi_1G_0\bar{x}+2\bar{x}_0^T\Psi_1B_0u_0
  +D_0^T\Psi_1D_0\big]dt,\\
    \label{eq33}
  &\mathbb{E}[\bar{x}^T(T)\bar{\Gamma}_0^TH_0\bar{\Gamma}_0\bar{x}(T)]-  \mathbb{E}[\bar{x}^T(0)\Psi_2(0)\bar{x}(0)]\cr
  =&\mathbb{E} \int_0^T\big[\bar{x}^T(\dot{\Psi}_2+\bar{A}^T\Psi_2+\Psi_2\bar{A})\bar{x}+2\bar{x}_0^T\bar{F}^T\Psi_2\bar{x}\big]dt,
  \end{align}
  and
   \begin{align}\label{eq34}
  &\mathbb{E}[\bar{x}^T(T)(-\bar{\Gamma}_0^TH_0)\bar{x}_0(T)]-  \mathbb{E}[\bar{x}^T(0)\Psi_3(0)\bar{x}_0(0)]\cr
  =&\mathbb{E} \int_0^T\big[\bar{x}^T(\dot{\Psi}_3+\bar{A}^T\Psi_3+\Psi_3{A}_0)\bar{x}_0+\bar{x}^T\Psi_3B_0u_0+\bar{x}_0^T\bar{F}^T\Psi_3\bar{x}_0\big]dt.
  \end{align}
  Note that by (\ref{eq3b}) and (\ref{eq18a}) we have
  $$d(x_0-\bar{x}_0)=[(A_0+B_0P_0)(x_0-\bar{x}_0)+G_0(x^{(N)}-\bar{x})]dt,$$
  which 
  with (\ref{eq91}) gives
  $\sup_{0\leq t\leq T}\mathbb{E}\|x_0-\bar{x}_0\|^2=O(1/N).$
From this and (\ref{eq32b})-(\ref{eq34}), one can obtain
  \begin{align}\label{eq97}
    \bar{J}_0(u_0,\hat{u})=& \mathbb{E}[\bar{x}_0^T(0)\Psi_1(0)\bar{x}_0(0)+\bar{x}^T(0)\Psi_2(0)\bar{x}(0)+2\bar{x}^T(0)\Psi_3(0)\bar{x}_0(0)]]\cr
    &+\mathbb{E}\int_0^T\Big[\bar{x}^T_0\Psi_1B_0R_0^{-1}B^T_0\Psi_1\bar{x}_0+\bar{x}^T\Psi_3B_0R_0^{-1}B^T_0\Psi_3\bar{x}+2\bar{x}^T\Psi_3B_0R_0^{-1}B^T_0\Psi_1\bar{x}_0 \cr
    &+2(\bar{x}_0^T\Psi_1+\bar{x}^T\Psi_3)B_0u_0+u_0^TR_0u_0+D_0^T\Psi_1D_0 \Big]dt\cr
        =& \mathbb{E}[\xi_0^T\Psi_1(0)\xi_0+\bar{\xi}^T\Psi_2(0)\bar{\xi}+2\bar{\xi}^T\Psi_3(0)\xi_0]\cr
        &+\mathbb{E}\int_0^T\Big[\big\|u_0+R_0^{-1}B^T_0\Psi_1\bar{x}_0
        +R_0^{-1}B^T_0\Psi_3\bar{x}\big\|^2_{R_0}+D_0^T\Psi_1D_0\Big]dt\cr
       \geq &\mathbb{E}[\xi_0^T\Psi_1(0)\xi_0+\bar{\xi}^T\Psi_2(0)\bar{\xi}+2\bar{\xi}^T\Psi_3(0)\xi_0]+\mathbb{E}\int_0^T\big(D_0^T\Psi_1D_0\big)dt.
  \end{align}
By (\ref{eq36a}) and (\ref{eq97}), we obtain (\ref{eq78}) and
 \begin{equation}\label{eq41}
   J_0(\hat{u}_0,\hat{u})\leq  \bar{J}_0(u_0,\hat{u})+O(1/\sqrt{N}).
    \end{equation}
Besides, it follows from (\ref{eq18a}) that
\begin{align*}
   \bar{J}_0({u}_0,\hat{u})=&\mathbb{E}\int_{0}^{T}\big[\|x_0- \Gamma_0 \hat{x}^{(N)}+\Gamma_0(\hat{x}^{(N)}-\bar{x})\|_{Q_0}^{2}+\|u_0\|_{R_0}^{2}\big]dt\cr
   &+\mathbb{E}[\|x_0(T)- \bar{\Gamma}_0\hat{x}^{(N)}(T)+\bar{\Gamma}_0(\hat{x}^{(N)}(T)-\bar{x}(T))\|_{H_0}^{2}]\cr
\leq&{J}_0({u}_0,\hat{u})+C\sup_{0\leq t \leq T}\mathbb{E}\Big[2\big(\|x_0- \Gamma_0\hat{x}^{(N)}\|^2\|\hat{x}^{(N)}-\bar{x}\|^2\big)^{1/2}+\|\hat{x}^{(N)}-\bar{x}\|^{2}\Big]\cr
\leq &{J}_0({u}_0,\hat{u})+O(1/\sqrt{N}).
\end{align*}
From this and (\ref{eq41}), we have
$J_0(\hat{u}_0,\hat{u})\leq {J}_0({u}_0,\hat{u})+O(1/\sqrt{N}).$
\hfill{$\Box$}

\begin{IEEEbiography}
{Bing-Chang Wang}
(SM'19) received the Ph.D. degree in System Theory from Academy of Mathematics and Systems Science, Chinese Academy of Sciences, Beijing, China, in 2011. From September 2011 to August 2012, he was with Department of Electrical and Computer Engineering, University of Alberta, Canada, as a Postdoctoral Fellow. From September 2012 to September 2013, he was with School of Electrical Engineering and Computer Science, University of Newcastle, Australia, as a Research Academic. From October 2013, he has been with School of Control Science and Engineering, Shandong University, China, and now is a Professor. He held visiting appointments as a Research Associate with Carleton University, Canada, from November 2014 to May 2015, and with the Hong Kong Polytechnic University from November 2016 to January 2017. His current research interests include mean field games, stochastic control, multiagent systems and reinforcement learning.
\end{IEEEbiography}

\begin{IEEEbiography}
	{Juanjuan Xu}  received the B.E. degree in mathematics from Qufu Normal University, Jining, China, in 2006, the M.E. degree in mathematics from Shandong University, Shandong, China, in 2009, and
	the Ph.D. degree in control science and engineering in
	2013 from Shandong University.
		
	She is currently a Qilu Professor with Shandong University. Her research interests include distributed consensus, optimal control,
	game theory, stochastic systems, and time-delay systems.
\end{IEEEbiography}

\begin{IEEEbiography}
{Huanshui Zhang}
(SM'06) received the B.S. degree in mathematics from Qufu Normal University, Shandong, China, in 1986, the M.Sc. degree in control theory from Heilongjiang University, Harbin, China, in 1991, and the Ph.D. degree in control theory from Northeastern University, Shenyang, China, in 1997.

He was a Postdoctoral Fellow at Nanyang Technological University, Singapore, from 1998 to 2001 and Research Fellow at Hong Kong Polytechnic University, Hong Kong, China, from 2001 to 2003. He currently holds a Professorship at Shandong University of Science and Technology, Qingdao, China. He was a Professor with the Harbin Institute of Technology, Shenzhen, China, from 2003 to 2006 and a Taishan Distinguished Professor and Changjiang Distinguished Professor with Shandong University, Jinan, China, from 2006 to 2019. He also held visiting appointments as a Research Scientist and Fellow with Nanyang Technological University, Curtin University of Technology, and Hong Kong City University from 2003 to 2006. His interests include optimal LQ control, decentralized control, time-delay systems, stochastic systems, Stackelberg game systems, and networked control systems. He was an Associate Editor for IEEE Transactions on Automatic Control and IEEE Transactions on Circuits and Systems I: Regular Papers.
\end{IEEEbiography}

\begin{IEEEbiography}
{Yong Liang} received the B.E. degree in building
electricity and intelligence from Qingdao University of Technology, Qingdao, China, in 2017, the
M.E. degree in control science and engineering
from Shandong University, Shandong, China, in
2020, the Ph.D. degree in control science and
control engineering from Shandong University,
Jinan, China, in 2024. He is currently working in
the School of Information Science and Engineering, Shandong Normal University, Jinan, China. 
His current research interests include mean
field games, social control, and networked control
systems.
\end{IEEEbiography}

\end{document}